\documentclass[showpacs,preprintnumbers,amsmath,amssymb,12pt]{article}
\usepackage{graphicx,amsfonts,amssymb,amsthm,amsmath,graphics,epsfig,pstricks,fancyhdr,fancybox}
\usepackage{dcolumn}
\usepackage{bm}

\newcommand{\refeq}[1]{~(\ref{#1})}

\newcommand{\totint}[2]{\int_{-\infty}^{+\infty}{#1}\,d{#2}}
\newcommand{\halfint}[2]{\int_{0}^{+\infty}{#1}\,d{#2}}

\newcommand{\Zeroint}[3]{\int_{{#1}\neq0}{#2}{#3}\,d{#1}}
\newcommand{\arcsinh}{\mathrm{arcsinh}\,}
\newcommand{\Real}{\mathbb{R}}
\newcommand{\rv}{\textit{rv}}
\newcommand{\id}{\textit{id}}
\newcommand{\sd}{\textit{sd}}
\newcommand{\ac}{\textit{ac}}
\newcommand{\pd}{pseudo--differential}
\newcommand{\LS}{\textit{L-S}}
\newcommand{\LKh}{L\évy--Khintchin}
\newcommand{\PDF}{\textit{pdf}}
\newcommand{\CHF}{\textit{chf}}
\newcommand{\LCH}{\textit{lch}}
\newcommand{\WF}{\textit{wf}}
\newcommand{\FT}{\textit{FT}}
\newcommand{\VG}{Variance--Gamma}
\newcommand{\Pqo}{\mathbf{P}\hbox{-a.s.}}
\newcommand{\eqd}{\stackrel{d}{=}}

\newcommand{\expect}[1]{\mathbf{E}\left[{#1}\right]}
\newcommand{\var}[1]{\mathbf{V}[{#1}]}
\newcommand{\rqm}{relativistic \textit{qm}}
\newcommand{\Rqm}{Relativistic \textit{qm}}
\newcommand{\pdf}{f}
\newcommand{\chf}{\varphi}
\newcommand{\lch}{\eta}
\newcommand{\scale}{a}
\newcommand{\tscale}{\tau}
\newcommand{\freq}{\omega}
\newcommand{\vel}{c}
\newcommand{\param}{\lambda}
\newcommand{\Nvar}{\scale}
\newcommand{\inNvar}{\sigma}
\newcommand{\prpdf}{p}
\newcommand{\prtrpdf}[2]{p({#1}|\,{#2})}
\newcommand{\prchf}{\phi}
\newcommand{\trpdf}{q}
\newcommand{\trchf}{\chi}
\newcommand{\inpdf}{f_0}
\newcommand{\inchf}{\varphi_0}
\newcommand{\inlch}{\eta_0}
\newcommand{\inwf}{\wf_0}
\newcommand{\inhatwf}{\hatwf_0}
\newcommand{\inscale}{b}
\newcommand{\complaw}{\mathfrak{H}}
\newcommand{\comppdf}{h}
\newcommand{\compchf}{\vartheta}
\newcommand{\complch}{\zeta}
\newcommand{\backpdf}{h_0}
\newcommand{\backchf}{\vartheta_0}
\newcommand{\backlch}{\zeta_0}
\newcommand{\wf}{\psi}
\newcommand{\hatwf}{\hat{\psi}}

\newcommand{\prop}{g}
\newcommand{\propchf}{\gamma}
\newcommand{\Ltriple}{\mathcal{L}}
\newcommand{\drift}{\alpha}
\newcommand{\diff}{\beta}
\newcommand{\Lpdf}{\ell}
\newcommand{\Lmeas}{\nu}
\newcommand{\testf}{v}
\newcommand{\gen}{A}

\newcommand{\law}{\mathfrak{F}}
\newcommand{\norm}{\mathfrak{N}}
\newcommand{\Poiss}{\mathfrak{P}}
\newcommand{\degen}{\mathfrak{D}}
\newcommand{\cauchy}{\mathfrak{C}}
\newcommand{\lapl}{\mathfrak{L}}
\newcommand{\unif}{\mathfrak{U}}
\newcommand{\stable}{\mathfrak{S}}
\newcommand{\vg}{\mathfrak{VG}}
\newcommand{\relat}{\mathfrak{R}}
\newcommand{\stud}{\mathfrak{T}}

\begin{document}
\thispagestyle{empty}

\title{\Huge \textbf{L\évy--Schr\"odinger wave packets}}
\author{\textsc{Nicola Cufaro Petroni}\\
Dipartimento di Matematica and \textsl{TIRES}, Universit\`a di Bari\\
\textsl{INFN} Sezione di Bari\\
via E. Orabona 4, 70125 Bari, Italy\\
\textit{email: cufaro@ba.infn.it}}

\date{}

\maketitle

\begin{abstract}
\noindent We analyze the time--dependent solutions of the pseudo--differential
L\évy--Schr\"o\-din\-ger wave equation in the free case, and we compare them with the associated
L\évy processes. We list the principal laws used to describe the time evolutions of both the L\évy
process densities, and the L\évy--Schr\"odinger wave packets. To have self--adjoint generators and
unitary evolutions we will consider only absolutely continuous, infinitely divisible L\évy noises
with laws symmetric under change of  sign of the independent variable. We then show several
examples of the characteristic behavior of the L\évy--Schr\"odinger wave packets, and in
particular of the bi-modality arising in their evolutions: a feature at variance with the typical
diffusive uni--modality of both the L\évy process densities, and the usual Schr\"odinger wave
functions.
\end{abstract}




\section{Introduction and notations}\label{notation}

In a recent paper~\cite{cufaro09} it has been shown how to extend the well known relation between
the Wiener process and the Schr\"odinger equation~\cite{nelson,feynman,bohmvigier,guerra} to other
suitable L\évy process. This idea -- discussed elsewhere only in the stable
case~\cite{garbaczewski1,laskin} -- leads to a \LS\ (L\évy--Schr\"odinger) equation containing
additional integral terms which take into account the possible jumping part of the background
noise. In fact, the infinitesimal generator of the Brownian semigroup (the Laplacian) being
substituted by the more general generator of a L\'evy semigroup, we get an integro-differential
operator with both a continuous (differential and Gaussian) and a jumping (integral, non Gaussian)
part. These ideas have already been discussed in the framework of \emph{stochastic
mechanics}~\cite{nelson,guerra} and are considered as a model for systems more general than just
the usual quantum mechanics: a true \emph{dynamical theory of L\évy processes} that can be applied
to several physical problems~\cite{applications}. The aim of this paper is now to show a number of
explicit examples of wave packets solutions of \LS\ free equations.

In recent years we have witnessed a considerable growth of interest
in non Gaussian stochastic processes -- and in particular into
L\'evy processes -- from statistical mechanics to mathematical
finance. In the physical field, however, the research scope is
presently rather confined to the stable processes and to the
corresponding fractional
calculus~\cite{garbaczewski1,laskin,mainardi}, while in the
financial domain a vastly more general type of processes is at
present in use. Here we suggest that a L\évy stochastic mechanics
should be considered as a dynamical theory of the entire gamut of
the \emph{infinitely divisible} processes with time reversal
invariance, and that the horizon of its applications should be
widened even to cases different from the quantum systems.

This approach has several advantages: first of all the use of general infinitely divisible
processes lends the possibility of having realistic, finite variances. Second, the presence of a
Gaussian component and the wide spectrum of decay velocities of the increment densities will give
the possibility of having models with differences from the usual Brownian (and usual quantum
mechanical, Schr\"odinger) case as small as we want. Last but not least, there are examples of non
stable L\'evy processes which are connected with the simplest form of the quantum,
\emph{relativistic} Schr\"odinger equation: an important link that was missing in the original
Nelson model. This final remark, on the other hand, shows that the present inquiry is not only
justified by the a desire of formal generalization, but is required by the need to attain
physically meaningful cases that otherwise would not be contemplated in the narrower precinct of
the stable laws.

In this paper we will show practical examples for the behavior of the evolving wave packet
solutions of particular kinds of (non Wiener) \LS\ equations, and we will put in evidence their
characteristics. In particular the \emph{bi--modality} arising in many of these these evolutions
which has a correspondence neither in the the process diffusions, nor in the usual Schr\"odinger
wave functions: an effect which has already been observed only in confined L\évy
flights~\cite{chechkin}. This is coherent with the usual stochastic mechanics scheme, in so far as
in this theory the Schr\"odinger equation is recovered by introducing a kind of interaction
modeled by means of a \emph{quantum potential}~\cite{nelson,guerra}. In the following exposition
laws and processes will always be one dimensional. An extensive analysis of the topics discussed
in this first chapter is available in the two monographs~\cite{sato} and~\cite{applebaum}, while a
short introduction can be found in~\cite{cufaro08}.

In the present paper the law of a \rv\ (random variable) $X$ is characterized either by its \PDF\
(probability density function) $\pdf$, when -- as it is generally supposed -- the law is \ac\
(absolutely continuous), or by its \CHF\ (characteristic function) $\chf$ with the usual
reciprocity relations
\begin{equation}\label{reciprocity}
    \chf(u)=\totint{\pdf(x)e^{iux}}{x},\qquad\quad
    \pdf(x)=\frac{1}{2\pi}\totint{\chf(x)e^{-iux}}{u}.
\end{equation}
When the laws are not \ac\ we sometimes will use the \emph{Dirac
delta} notation: the symbol $\delta_{x_0}(x)=\delta(x-x_0)$ will
then represent a law degenerated in $x_0$ and only formally it will
act as a \PDF. The symbol $\delta(x)$ will also be used instead of
$\delta_0(x)$. In order to have background noises with generators
self--adjoint in $L^2$ -- an essential requirement for our purposes
-- we will consider only \emph{symmetric} laws, namely we will
require
\begin{equation*}
    \pdf(-x)=\pdf(x)\,,\qquad\quad\chf(-u)=\chf(u)
\end{equation*}
so that the \CHF\ $\chf$ will also be real. This also means that, when it exists,  the expectation
vanishes $\expect{X}=0$, namely the law is also \emph{centered}. See the Appendix~\ref{types} for
further details about our notations.

Since we will restrict our analysis to background noises driven by
L\évy processes, we will be interested almost exclusively in \id\
(infinitely divisible\footnote{A law $\varphi$ is said to be \id\ if
for every $n$ it exists a \CHF\ $\varphi_n$ such that
$\varphi=\varphi_n^n$; on the other hand $\chf$ is said to be stable
when for every $c>0$ it is always possible to find $a>0$ and
$b\in\mathbf{R}$ such that $e^{ibu}\varphi(au)=[\varphi(u)]^c$.
Every stable law is also \id. See also Appendix~\ref{types} for
further details about stable laws.}, for details
see~\cite{sato,applebaum,cufaro08}) laws with a L\évy triplet
$\Ltriple=(\drift,\diff,\Lmeas)$. Here our L\évy measures $\Lmeas$
will always be supposed to have a density:
$\Lmeas(dy)=\Lpdf(y)\,dy$; when this does not happen we will often
use the Dirac delta notation. As a consequence the L\évy triplet
will be rather specified as $\Ltriple=(\drift,\diff,\Lpdf)$. The
\LCH\ (logarithmic characteristic) of our \id\ laws $\lch=\ln\chf$,
with $\chf=e^\lch$, will then satisfy the \LKh\ formula
\begin{equation}\label{LKh}
    \lch(u)=i\drift
    u-\frac{1}{2}\,\diff^2u^2+\Zeroint{y}{\left[e^{iuy}-1-iuy\,I_D(y)\right]}{\,\Lpdf(y)}
\end{equation}
where $D=\{y\,:\, |y|<1\}$. The prescription of the integral around the origin is essential only
when -- as usually may happen -- the L\évy measure shows a singularity in $y=0$. When the law is
dimensionless (see Appendix~\ref{types}) then also $\drift,\,\diff,\,\Lpdf$ and $y$ are so; on the
other hand, if the law has the dimensions of a length, then $\drift,\diff,y$ are lengths, while
$\Lpdf$ is the reciprocal of a length. In particular when the law is symmetric we have
\begin{equation*}
    \drift=0\,,\qquad\quad \Lpdf(-x)=\Lpdf(x)
\end{equation*}
so that the \LKh\ formula will be reduced to the symmetric real expression
\begin{equation}\label{sLKh}
    \lch(u)=-\frac{1}{2}\,\diff^2u^2+\Zeroint{y}{(\cos uy-1)}{\,\Lpdf(y)}
\end{equation}
and hence the \CHF\ $\chf$ will not only be real, but also non negative: $\chf(u)\geq0$.


The Markov processes dealt with in this paper are stationary,
independent increments processes and are then defined by means of
the \CHF\ $\chf^{\Delta t/\tau}$ of their $\Delta t$--increments,
where $\tau$ is a dimensional, time scale parameter. Here too we can
introduce a dimensionless formulation through the coordinate
$t/\tau$, and to simplify the notation we can continue to use the
same symbol $t$ for this dimensionless time. In this case the
stationary \CHF\ will be reduced to $\chf^{\Delta t}$, and the
dimensional formulation will be recovered by simple substitution of
$t/\tau$ to $t$. A stochastically continuous process with stationary
and independent increments is called a \emph{L\évy process} when
$X(0)=0,\;\Pqo$, but this paper will mostly be about the same kind
of processes for arbitrary initial conditions $X(0)=X_0,\;\Pqo$ with
law $\inpdf(x)$ and $\inchf(u)=e^{\inlch(u)}$. All these processes,
independently from their initial conditions, will share both the
same differential equations (whether \emph{SDE}'s, or \emph{PDE}'s)
and the same transition \PDF's
\begin{equation*}
    \pdf_{X(t)}\left(x\,|\,X(s)=y\right)\,\,=\,\prtrpdf{x,t}{y,s}.
\end{equation*}
To avoid confusion we will then adopt different notations for their
respective marginal \PDF's: for a L\évy process (namely with $X_0=0$
initial condition) we will write
\begin{equation*}
    \pdf_{X(t)}(x)=\trpdf(x,t),\qquad\quad\chf_{X(t)}(u)=\trchf(u,t)
\end{equation*}
with $\trpdf(x,0)=\delta(x)$ and $\trchf(x,0)=1$, while for the general stationary and independent
increments process (with arbitrary initial condition $X_0$) we will write
\begin{equation*}
    \pdf_{X(t)}(x)=\prpdf(x,t),\qquad\quad\chf_{X(t)}(u)=\prchf(u,t)
\end{equation*}
with $\prpdf(x,0)=\inpdf(x)$ and $\prchf(x,0)=\inchf(x)$. It is then
easy to show that
\begin{equation}\label{prtrpdf}
  \prtrpdf{x,t}{y,s}=\trpdf(x-y,t-s).
\end{equation}

The infinitesimal generator $\gen=\lch(\partial)$ (here $\partial$
stands for the derivation with respect to the variable of a test
function $\testf$) of the semigroup of a L\évy process will be a
\pd\ operator with symbol $\lch$~\cite{cufaro09,applebaum}, namely
from\refeq{LKh}
\begin{eqnarray}\label{pdgen}
    [\gen\testf](x)&=&[\lch(\partial)\testf](x)=\frac{1}{\sqrt{2\pi}}\totint{e^{iux}\lch(u)\hat{\testf}(u)}{u}\nonumber\\
    &=&\drift\,\partial_x\testf(x)+\frac{\diff^2}{2}\,\partial_x^2\testf(x)\\
    &&
    \qquad\qquad\qquad\quad+\Zeroint{y}{\left[\testf(x+y)-\testf(x)-yI_D(y)\,\partial_x\testf(x)\right]}{\,\Lpdf(y)}\nonumber
\end{eqnarray}
where $\hat \testf$ denotes the \FT\ (Fourier transform) of the test
function $\testf$ with the reciprocity relations:
\begin{equation*}
    \hat\testf(u)=\frac{1}{\sqrt{2\pi}}\totint{\testf(x)e^{-iux}}{x},\qquad\quad
    \testf(x)=\frac{1}{\sqrt{2\pi}}\totint{\hat\testf(u)e^{iux}}{u}
\end{equation*}
The generator $\gen$ will be self--adjoint in $L^2(\mathbb{R},dx)$
when the law is symmetric, and in this case it reduces to
\begin{equation}\label{pdsymmgen}
    [\gen\testf](x)=\frac{\diff^2}{2}\,\partial_x^2\testf(x)+\Zeroint{y}{\left[\testf(x+y)-\testf(x)\right]}{\,\Lpdf(y)}
\end{equation}
so that it is determined by the two essential elements of our L\évy
triplet, namely $\diff$ and $\Lpdf$.
\begin{table}
  \centering
  \begin{tabular}{|c||c|c||c|c||c|c|}
  \hline
  \textit{law} & $\pdf$ & $\chf$ & $\diff$ & $\ell$ & $\mathbf{E}$ & $\mathbf{V}$ \\ \hline \hline
  & & & & & &\\
  $\degen$ & $\delta(x)$ & $1$ & $0$ & $0$ & 0& 0\\
  & & & & & &\\
  $\norm$ & $\frac{e^{-x^2/2}}{\sqrt{2\pi}}$ & $e^{-u^2/2}$ & $1$ & $0$ & 0& 1\\
  & & & & & &\\
  $\cauchy$ & $\frac{1}{\pi}\,\frac{1}{1+x^2}$ & $e^{-|u|}$ & $0$ & $\frac{1}{\pi x^2}$  & --&$+\infty$\\
  & & & & & &\\
  $\lapl$ & $\frac{e^{-|x|}}{2}$ & $\frac{1}{1+u^2}$ & $0$ & $\frac{e^{-|x|}}{|x|}$  &0 &2\\
  & & & & & &\\
  $\unif$ & $\frac{\Theta(x+1)-\Theta(x-1)}{2}$ & $\frac{\sin u}{u}$ & -- & --  & 0&$\frac{1}{3}$\\
  & & & & & &\\
  $\degen_1$ & $\frac{\delta_1(x)+\delta_{-1}(x)}{2}$ & $\cos u$ & -- & --  & 0&$1$\\
  & & & & & &\\
  \hline
  \end{tabular}
  \caption{List of the essential properties of a few basic, dimensionless
laws discussed in this paper: degenerate (Dirac) $\degen$, normal (Gauss) $\norm$, Cauchy
$\cauchy$, Laplace $\lapl$, uniform $\unif$, and doubly degenerate in $+1,-1$ (symmetric
Bernoulli) $\degen_1$.}\label{basiclaws}
\end{table}
Given the process stationarity, in a dimensionless formulation the
transition law degenerate in $x=0$ at $t=0$ will have as \CHF\
$\trchf=\chf^t=e^{t\lch}$ and as \PDF
\begin{equation}\label{trpdf}
    \trpdf(x,t)=\frac{1}{2\pi}\totint{\trchf(u,t)e^{-iux}}{u}=
               \frac{1}{2\pi}\totint{\chf(u)^te^{-iux}}{u}.
\end{equation}
This transition law plays an important role in the evolution of an
arbitrary initial law $\inpdf, \;\inchf$: the process \CHF\ will
indeed be now $\prchf(u,t)=\trchf(u,t)\inchf(u)$, and the
corresponding \PDF\ will be calculated from
\begin{equation*}
    \prpdf(x,t)=[\trpdf(t)*\inpdf](x)=\frac{1}{2\pi}\totint{\prchf(u,t)e^{-iux}}{u}
\end{equation*}
namely either as a convolution of the transition and the initial \PDF's, or by inverting the \CHF\
$\prchf$ of the process. This \PDF\ will also be a solution of the evolution \pd\
equation~\cite{cufaro09,applebaum}
\begin{equation}\label{pdfeq}
    \partial_t\prpdf=\eta(\partial)\prpdf,\qquad\quad\prpdf(x,0)=\inpdf(x)
\end{equation}
which from\refeq{pdgen} takes the integro--differential form
\begin{eqnarray*}
    \partial_t\prpdf(x,t)&=&\drift\,\partial_x\prpdf(x,t)+\frac{\diff^2}{2}\,\partial_x^2\prpdf(x,t)\\
    &&\qquad\qquad +\Zeroint{y}{\left[\prpdf(x+y,t)-\prpdf(x,t)-yI_D(y)\,\partial_x\prpdf(x,t)\right]}{\,\Lpdf(y)}
\end{eqnarray*}
and for a centered, symmetric noise from\refeq{pdsymmgen} reduces to
\begin{equation}\label{symmFP}
    \partial_t\prpdf(x,t)=\frac{\diff^2}{2}\,\partial_x^2\prpdf(x,t)+\Zeroint{y}{\left[\prpdf(x+y,t)-\prpdf(x,t)\right]}{\,\Lpdf(y)}\,.
\end{equation}
We finally remember that, since\refeq{pdfeq} and\refeq{symmFP} are
given in terms of process \PDF's, this equations are supposed to
hold only for \ac\ processes. We are then required to point out
which L\évy processes have densities. To answer -- at least
partially -- this question we then recall that~\cite{sato} any
non--degenerate, \sd\ (self--decomposable\footnote{A law $\chf(u)$
is \sd\ when for every $a\in(0,1)$ we can always find another \CHF\
$\chf_a(u)$ such that $\chf(u)=\chf(au)\chf_a(u)$. Every stable law
is also \sd; every \sd\ law is also \id.}, for details
see~\cite{sato,applebaum,cufaro08}) distribution is \ac. On the
other hand such a property also extends to the corresponding
processes for every $t$. In fact~\cite{sato} if $X(t)$ is a \sd\
process also its \PDF\ at every $t$ is \sd, and hence $X(t)$ is \ac\
for every $t$. As a consequence we can always explicitly write down
the evolution equations\refeq{symmFP} in terms of the process \PDF's
at least for the \sd\ case. We remark, however, that there are also
non \sd\ processes which are \ac: the \ac\ compound Poisson
processes of Appendix~\ref{poisson} are an example in point.

We listed in the Table~\ref{basiclaws} the properties of a few basic, symmetric, dimensionless
laws: degenerate (Dirac) $\degen$, normal (Gauss) $\norm$, Cauchy $\cauchy$, Laplace $\lapl$,
uniform $\unif$, and doubly degenerate in $+1,-1$ (symmetric Bernoulli) $\degen_1$. The uniform
law \PDF\ is given by means of the Heaviside functions $\Theta(x)$.
These laws are also relevant particular cases of the families that we will introduce in the
Section~\ref{SDLaws}.
Remark that in the Table~\ref{basiclaws} there is no value for the expectation of $\cauchy$
because it does not exist ($\cauchy$ is centered on the median), and no values for the L\évy
triplet of $\unif$ and $\degen_1$ since these are not \id\ laws. Moreover in general our laws are
not necessarily standard. The form of the simplest generators corresponding to our L\évy processes
is finally shown in the Table~\ref{basiclawsgen}.
\begin{table}
  \centering
  \begin{tabular}{|c||c|}
  \hline
  \textit{law} & $[\gen\testf](x)$ \\ \hline \hline
  & \\
  $\norm$ & $\frac{\partial_x^2 \testf(x)}{2}$ \\
  & \\
  $\cauchy$ & $\Zeroint{y}{\frac{\testf(x+y)-\testf(x)}{\pi y^2}}{}$ \\
  & \\
  $\lapl$ & $\Zeroint{y}{\frac{[\testf(x+y)-\testf(x)]e^{-|y|}}{|y|}}{}$ \\
  & \\
  \hline
  \end{tabular}
  \caption{List of the generators of the L\évy processes associated to some of
  the non degenerate, \id, dimensionless laws of Table~\ref{basiclaws}}\label{basiclawsgen}
\end{table}

The paper is organized as follows: in the Chapter~\ref{SDLaws} we recall the essential properties
of the law families of our interest; then in the Chapter~\ref{LS} the \LS\ equation is introduced
with its connections to the L\évy processes. Finally in Chapter~\ref{wpack} our examples are
elaborated and in Chapter~\ref{conclusions} the results are collected and discussed. A few
technical details are collected in the Appendices in order to avoid to excessively burden the
text.

\section{Families of \id\ laws}\label{SDLaws}

\begin{figure}
\begin{center}
\vspace{-1.0cm}
\includegraphics*[width=14.0cm]{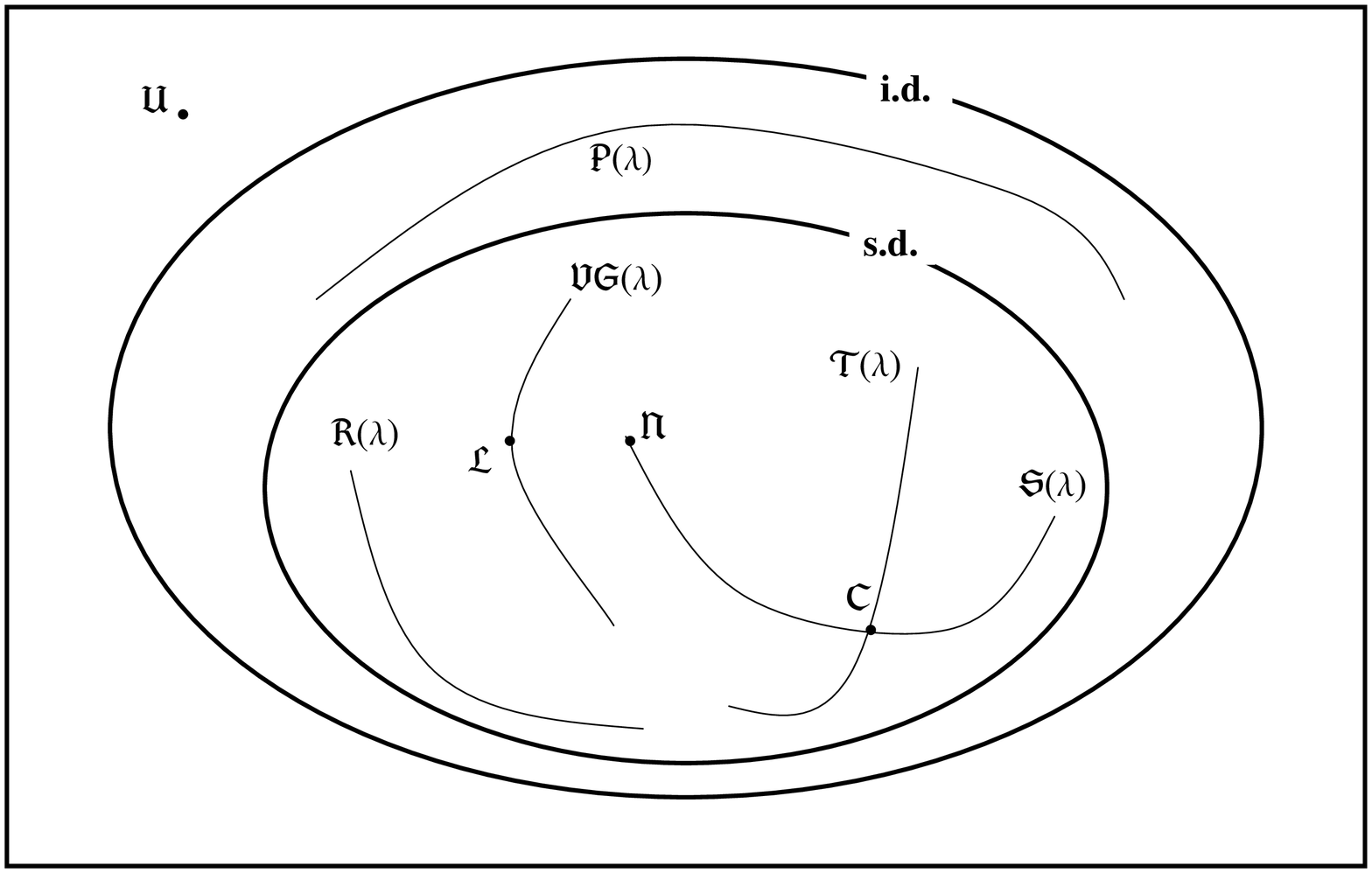}
\end{center}
\vspace{-1.5cm} \caption{Graphical synthesis of the relations among the families of laws discussed
in the Section~\ref{SDLaws}. The uniform $\unif$ is our unique example beyond the pale of the \id\
laws, while the laws of the (simple) Poisson family $\Poiss(\param)$ are \id\ but not \sd. Notable
cases ($\norm$, $\lapl$, $\cauchy$) within the \sd\ families are put in evidence; the Cauchy
$\cauchy$ law lies at the intersection of the stable $\stable(\param)$ and Student $\stud(\param)$
families.}\label{schema}
\end{figure}
\begin{figure}
\begin{center}
\includegraphics*[width=13.0cm]{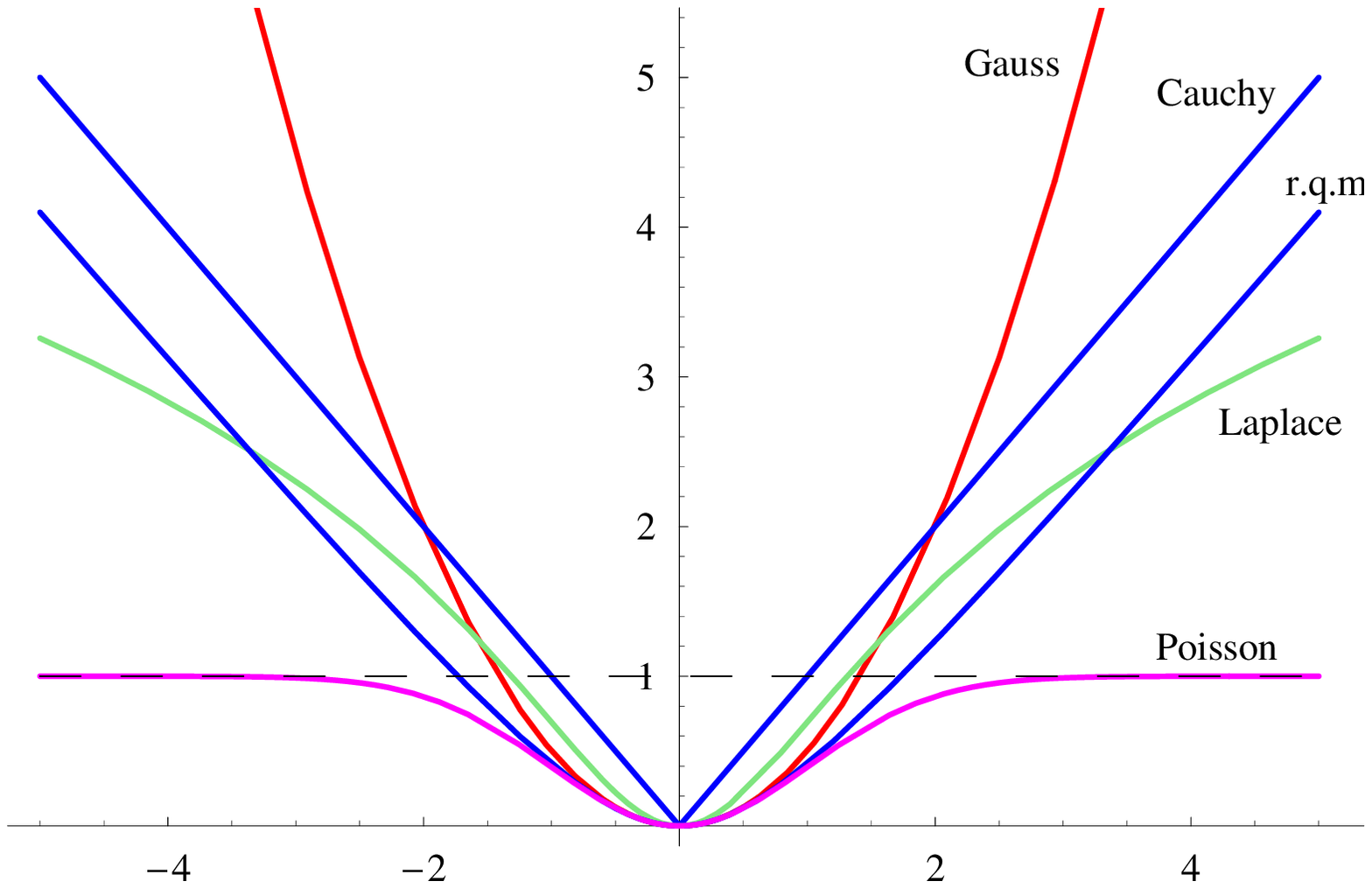}
\end{center}
\caption{The \LCH\ $-\eta(u)$ of some basic dimensionless laws from Table~\ref{families}, plus
that of a compound Poisson $\Poiss(\param,\norm)$ with normal component laws (see
Appendix~\ref{poisson}).}\label{fig01}
\end{figure}
\begin{table}
  \centering
  \begin{tabular}{|c||c|c||c|c||c|}
  \hline
  \textit{law} & $\pdf$ & $\chf$ & $\diff$ & $\ell$ & $\param>0$  \\ \hline \hline
  & & & & & \\
  $\stable(\param)$ & $\begin{array}{c}
                        \vspace{5pt} H_\param(|x|)   \\ \vspace{5pt}
                       \frac{1}{\pi}\,\frac{1}{1+x^2} \\
                       \frac{e^{-x^2/2}}{\sqrt{2\pi}}
                       \end{array}$
                     & $\begin{array}{c}
                        \vspace{5pt} e^{-|u|^\param/\param} \\ \vspace{5pt}
                       e^{-|u|} \\
                       e^{-u^2/2}
                       \end{array}$
                     & $\begin{array}{c}
                        \vspace{5pt} 0 \\ \vspace{5pt}
                       0 \\
                       1
                       \end{array}$
                     & $\begin{array}{c}
                        \vspace{5pt} \frac{|x|^{-1-\param}}{-2\param\Gamma(-\param)\cos(\param\pi/2)} \\ \vspace{5pt}
                       \frac{1}{\pi x^2} \\
                       0
                       \end{array}$
                     & $\begin{array}{c}
                        \vspace{5pt} <2  \\ \vspace{5pt}
                       1 \\
                       2
                       \end{array}$
                     \\
  & & & & & \\
  $\vg(\param)$ & $\begin{array}{c}
                       \vspace{5pt}  \frac{|x|^{\param-1/2}K_{\param-1/2}(|x|)}{2^{\param-1}\Gamma(\param)\sqrt{2\pi}} \\
                       \frac{e^{-|x|}}{2}
                       \end{array}$
                 & $\begin{array}{c}
                       \vspace{5pt} \left(\frac{1}{1+u^2}\right)^\param \\
                        \frac{1}{1+u^2}
                       \end{array}$
                 & $\begin{array}{c}
                        \vspace{5pt} 0  \\
                       0
                       \end{array}$
                 & $\begin{array}{c}
                        \vspace{5pt} \frac{\param\, e^{-|x|}}{|x|}  \\
                       \frac{e^{-|x|}}{|x|}
                       \end{array}$
                 & $\begin{array}{c}
                        \vspace{5pt} \ldots \\
                       1
                       \end{array}$
                 \\
  & & & & & \\
  $\stud(\param)$ & $\begin{array}{c}
                        \vspace{5pt} \frac{1}{B\left(\frac{1}{2},\frac{\param}{2}\right)}\left(\frac{1}{1+x^2}\right)^{\frac{\param+1}{2}}  \\
                       \frac{1}{\pi}\,\frac{1}{1+x^2}
                       \end{array}$
                   & $\begin{array}{c}
                        \vspace{5pt} \frac{2|u|^{\param/2}K_{\param/2}(|u|)}{2^{\param/2}\Gamma(\param/2)}  \\
                       e^{-|u|}
                       \end{array}$
                   & $\begin{array}{c}
                        \vspace{5pt} 0  \\
                       0
                       \end{array}$
                   & $\begin{array}{c}
                        \vspace{5pt} \ldots  \\
                       \frac{1}{\pi x^2}
                       \end{array}$
                   & $\begin{array}{c}
                        \vspace{5pt} \ldots \\
                       1
                       \end{array}$
                   \\
  & & & & & \\
    $\relat(\param)$ & $\frac{\param e^\param K_1\left(\sqrt{\param^2+x^2}\right)}{\pi\sqrt{\param^2+x^2}}$
                     & $e^{\,\param(1-\sqrt{1+u^2})}$ & $0$ & $\frac{\param K_1(|x|)}{\pi |x|} $  & $\ldots$ \\
  & & & & &\\
  \hline
  \end{tabular}
  \caption{Properties of our principal families of \sd, dimensionless laws: the stable
$\stable(\param)$, the \VG\ $\vg(\param)$, the Student
$\stud(\param)$ and the \rqm\ (quantum mechanics)
$\relat(\param)$.}\label{families}
\end{table}
We will introduce here the principal families of \id\ laws considered in this paper. For a
graphical synthesis of the relations among them see Figure~\ref{schema}. Please remark that this
synthesis is particularly simple because we limit ourselves here to \emph{dimensionless} laws (see
Appendix~\ref{types}): this produces one-parameter families that can be easily represented in our
scheme. In the Table~\ref{families} are then listed the properties of the principal families of
dimensionless, \sd\ laws that will be discussed. The $"\ldots"$ symbol in this table means either
that we do not have an elementary formulation for the entry, or that there are no particular
values of $\param$ to be put in evidence. $K_\nu$, $B$ and $\Gamma$ respectively are the modified
Bessel functions of the second kind, and the Euler Beta and Gamma functions, while $H_\param$
stands for the Fox $H$-functions representing the \PDF\ of stable laws~\cite{schneider}. From
Table~\ref{basiclaws} and Table~\ref{families} we can on the other hand immediately see that
$\stable(1)=\stud(1)=\cauchy$, $\;\stable(2)=\norm$ and $\vg(1)=\lapl$, as also put in evidence in
the Figure~\ref{schema}. The behaviors of a few \LCH's of \id\ laws are finally displayed and
compared in the Figure~\ref{fig01}: it could be seen there that all the \LCH's of the \sd\ laws
considered in this paper diverge at infinity with velocities ranging from $u^2$ to $\log u$, while
the unique not diverging \LCH\ characterizes one of our non \sd\ examples: the compound Poisson
$\Poiss(\param,\norm)$ with normal component laws. This also gives an intuitive idea of how much
the behavior of a law -- in so far as we are concerned, for instance, with its jumping properties
-- differs from that of the  $\Pqo$ continuous Gaussian case.

\subsection{The stable laws $\stable(\param)$}

This is the more widely studied family of \id\ laws, albeit among
them only the normal $\stable(2)=\norm$ enjoys a finite variance.
But for the $\norm$, the $\cauchy$ and precious few other cases the
\PDF's of the stable laws exist only in the form of Fox
$H$-functions~\cite{schneider}.
To see in what sense these laws are \emph{stable} we must for a moment reintroduce the dimensional
parameter $a$: we then have a larger family $\stable_\scale(\param)$ with two parameters,
$0<\param\leq2$ and $\scale>0$, and
\begin{equation}\label{stable}
    \chf(u)=e^{-a^\param|u|^\param/\param}.
\end{equation}
Now, for a given fixed $\param$, the family $\stable_\scale(\param)$
with $a>0$ is closed under convolution, as can be easily seen
from\refeq{stable}. For instance the families of the normal
$\norm_\scale=\norm(\scale^2)$ and Cauchy $\cauchy_\scale$ laws are
closed under convolution since
$\norm(a_1^2)\!*\norm(a_2^2)=\norm(a_1^2+a_2^2)$ and
$\cauchy_{\scale_1}*\cauchy_{\scale_2}=\cauchy_{\scale_1+\scale_2}$.
Stability however means more: the families $\stable_\scale(\param)$
for a given $\param$ are \emph{types} of laws, in the sense that a
law of the family differs from another just by a re-scaling
(centering is not necessary here because our laws already are
centered; for details see Appendix~\ref{types}), the parameter $a$
being indeed nothing else than a space scale parameter. This has far
reaching consequences. In particular it is at the root of the well
known fact that the stable L\évy processes are \emph{self--similar}:
a property not extended to other, non stable L\évy
processes~\cite{cufaro08}. The generators of the stable L\évy
processes are
\begin{equation*}
    [\gen\testf](x)=\frac{-1}{2\param\Gamma(-\param)\cos\frac{\param\pi}{2}}\Zeroint{y}{\frac{\testf(x+y)-\testf(x)}{|y|^{1+\param}}}{}
    \qquad0<\param<2,\quad\param\neq1
\end{equation*}
while for $\param=1$ ($\cauchy$ law) and $\param=2$ ($\norm$ law)
they are listed in the Table~\ref{basiclawsgen}.

\subsection{The \VG\ laws $\vg(\param)$}

The \VG\ laws owe their name to the fact that they can be seen as
\emph{normal variance-mean mixtures}\footnote{A normal variance-mean
mixture, with mixing probability density $g$, is the law of a random
variable $Y$ of the form $Y=\alpha + \diff V+\sigma \sqrt{V}X$ where
$\alpha$ and $\diff$ are real numbers and $\sigma > 0$. The random
variables $X$ and $V$ are independent; $X$ is a normal standard, and
$V$ has a \PDF\ $g$ with support on the positive half-axis. The
conditional distribution of $Y$ given $V$ is then a normal
distribution with mean $\alpha + \diff V$ and variance $\sigma^2V$.
A normal variance-mean mixture can be thought of as the distribution
of a certain quantity in an inhomogeneous population consisting of
many different normally distributed sub--populations.} where the
mixing density is a \emph{gamma} distribution. It is apparent
moreover from the Table~\ref{families} that $\vg(\param)$ is closed
under convolution in the sense that
$\vg(\param_1)*\vg(\param_2)=\vg(\param_1+\param_2)$. That
notwithstanding, however, the \VG\ laws are not stable. To see that
let us reintroduce the dimensional scale parameter $\scale$ to have
the enlarged family $\vg_\scale(\param)$:
\begin{equation*}
  \chf(u)=\left(\frac{1}{1+\scale^2u^2}\right)^\param
\end{equation*}
Now every sub--family with a given, fixed $\scale$ is closed under convolution, but at variance
with the stable case the parameter describing the sub--family is $\param$, rather than $\scale$.
As a consequence the closed subfamilies do not constitute types of laws differing only by a
rescaling, and hence the laws are not stable. The \PDF's of the \VG\ laws can be given in
particular instances as finite combinations of elementary functions. By generalizing the quoted
example of the Laplace law $\vg(1)=\lapl$, when $\param=n+1$ with $n=0,1,\ldots$ we have for the
dimensionless \PDF's
\begin{equation*}
    \pdf(x)=\sum_{k=0}^n\binom{2n-k}{n}\frac{(2|x|)^ke^{-|x|}}{k!\,2^{2n+1}}=\frac{e^{-|x|}}{n!2^{n+1}}\,\theta_n(|x|)
\end{equation*}
where $\theta_n(x)$ are reverse Bessel polynomials~\cite{grosswald}.
All our dimensionless $\vg(\param)$ laws are endowed with
expectations (which vanish by symmetry) and finite variances
$2\param$. The generator of the corresponding L\évy process is
\begin{equation*}
    [\gen\testf](x)=\param\Zeroint{y}{\frac{\testf(x+y)-\testf(x)}{|y|}}{\,e^{-|y|}}\qquad\quad
    \param>0
\end{equation*}
which coincides with that of $\lapl$ (see Table~\ref{basiclawsgen})
for $\param=1$.

\subsection{The Student laws $\stud(\param)$}

But for the Cauchy $\cauchy$ case, the laws of the Student family (even enlarged by means of the
scale parameter $\scale$) are not stable, and $\stud(\param)$ itself is not closed under
convolution: convolutions of Student laws are not Student laws. As can be seen from
Table~\ref{families} the \VG\ and the Student families enjoy a sort of duality since their \PDF's
and \CHF's are essentially exchanged. This has been discussed at length in a few recent
papers~\cite{cufaro07,BergVignat,HeideLeonenko}. Remark that to put in evidence this
correspondence we have chosen the Student laws of $\stud(\param)$ without introducing the usual
parametric scaling $x^2/\param$ of its variable that would have put equal to $\param/(\param-2)$
all their variances for $\param>2$. In particular this means that for $\param\to+\infty$ we will
not get a standard $\norm$ law, as also shown in the Figure~\ref{schema}. The following remarks
are however virtually untouched by this choice. While the \PDF's and \CHF's of the Student laws
are known, differently from the \VG\ laws their L\évy measures and generators have not a known
general expression. However we can give them in particular instances. For example when
$\param=2n+1$ with $n=0,1,\ldots$ the \CHF\ becomes
\begin{equation*}
    \chf(u)=\sum_{k=0}^n\frac{n!(2n-k)!}{(2n)!(n-k)!k!}(2|u|)^ke^{-|u|}=\frac{n!2^ne^{-|u|}}{(2n)!}\,\theta_n(|u|).
\end{equation*}
where $\theta_n$ are again reverse Bessel
polynomials~\cite{grosswald}. Of course $\stud(1)=\cauchy$ is the
well known Cauchy (stable) case, while for $\stud(3)$ we have
\begin{equation*}
    \pdf(x)=\frac{2}{\pi}\left(\frac{1}{1+x^2}\right)^2,\qquad\quad
    \chf(u)=(1+|u|)e^{-|u|}.
\end{equation*}
and it can be shown~\cite{cufaro07} in this case that the L\évy
measure is
\begin{equation}\label{studentLF}
    \ell(x)=\frac{1-|x|(\sin|x|\,\mathrm{ci}\,|x|-\cos|x|\,\mathrm{si}\,|x|)}{\pi x^2}
\end{equation}
where the sine and the cosine integral functions are
\begin{equation*}
    \mathrm{si}\, x= -\int_x^{+\infty}\frac{\sin y}{y}\,dy\,,\qquad\quad\mathrm{ci}\, x= -\int_x^{+\infty}\frac{\cos
    y}{y}\,dy
\end{equation*}
The existence of the moments of the $\stud(\param)$ laws depends on
the value of the parameter $\param$: the $n^{\mathrm{th}}$ moment
exists if $n<\param$. In particular the expectation exists (and
vanishes) for $\param>1$, while the variance exists finite for
$\param>2$ and its value is $(\param-2)^{-1}$. The generator of the
L\évy process can finally be explicitly given for $\stud(3)$
from\refeq{studentLF}
\begin{equation*}
    [\gen\testf](x)=\param\Zeroint{y}{[\testf(x+y)-\testf(x)]}{\,\frac{1-|y|(\sin|y|\,\mathrm{ci}\,|y|-\cos|y|\,\mathrm{si}\,|y|)}{\pi
    y^2}}.
\end{equation*}

\subsection{The compound Poisson laws $\norm_{\inNvar}*\Poiss(\param,\complaw)$}\label{poisslaw}

The compound Poisson laws  $\complaw_0*\Poiss(\param,\complaw)$ are
not \sd, but they are nevertheless \id; they are also \ac\ when
$\complaw_0$ is \ac\ (for details and notations see
Appendix~\ref{poisson}). In the following examples we will take into
account the dimensional parameters $\scale$ of the component laws.
Consider now the case $\complaw_0=\norm_{\inNvar}$: then
$\Lpdf_0(x)=0$ and $\diff_0=\inNvar$ so that the L\évy triplet of
$\norm_{\inNvar}*\Poiss(\param,\complaw)$ is $
\Ltriple=\left(0\,,\,\inNvar\,,\,\param\,\comppdf\right)$ and the
generator is
\begin{equation*}
    [\gen\testf](x)=\frac{\inNvar^2}{2}\,\partial_x^2\testf(x)+\param\Zeroint{y}{[\testf(x+y)-\testf(x)]}{\comppdf(y)}\,.
\end{equation*}
When in particular also $\complaw=\norm_\Nvar$, then the L\évy triplet of
$\norm_{\inNvar}*\Poiss(\param,\norm_\Nvar)$ is
\begin{equation*}
    \Ltriple=\left(0\,,\,\inNvar\,,\,\param\,\frac{e^{-x^2/2\Nvar^2}}{\sqrt{2\pi\Nvar^2}}\right)
\end{equation*}
and we get a law with the following \PDF\ and \LCH
\begin{equation*}
    \pdf(x)=e^{-\param}\sum_{k=0}^{\infty}\frac{\param^k}{k!}\,\frac{e^{-x^2/2(k\Nvar^2+\inNvar^2)}}{\sqrt{2\pi(k\Nvar^2+\inNvar^2)}},
    \quad\qquad\lch(u)=\param(e^{-\Nvar^2u^2/2}-1)-\frac{\inNvar^2u^2}{2},
\end{equation*}
namely a Poisson mixture of centered normal laws
$\norm(k\Nvar^2+\inNvar^2)$. The self--adjoint generator then is
\begin{equation*}
    [\gen\testf](x)=\frac{\inNvar^2}{2}\,\partial_x^2\testf(x)+
    \param\totint{[\testf(x+y)-\testf(x)]\frac{e^{-y^2/2\Nvar^2}}{\sqrt{2\pi\Nvar^2}}}{y}\,.
\end{equation*}
and we could look at it as to a Poisson correction to the Wiener
generator, the relative weight of these two independent components
being ruled by the ratio between $\param$ and $\inNvar^2$.

As another example of \ac\ compound Poisson law let us suppose
instead that $\complaw_0=\norm_{\inNvar}$ again, but that
$\complaw=\degen_\scale$ (see Appendix~\ref{poisson}), so that the
L\évy triplet of $\norm_{\inNvar}*\Poiss(\param,\degen_\scale)$ now
is
\begin{equation*}
    \Ltriple=\left(0\,,\,\inNvar\,,\,\param\,\frac{\delta_1(x/\scale)+\delta_{-1}(x/\scale)}{2\scale}\right)
\end{equation*}
while its \PDF\ and \LCH\ are
\begin{eqnarray*}
  \pdf(x) &=& e^{-\param}\sum_{k=0}^{\infty}\frac{\param^k}{k!}\,\frac{1}{ 2^k}\sum_{j=0}^k\binom{k}{j}
                                \frac{e^{-[x-(k-2j)\scale]^2/2\inNvar^2}}{\sqrt{2\pi\inNvar^2}} \\
  \lch(u) &=&\param(\cos \scale u-1)-\frac{\inNvar^2u^2}{2}
\end{eqnarray*}
Here the law is again a mixture of normal laws
$\norm(n\scale,\inNvar^2),\; n=0,\pm1,\ldots$ which however are now
centered around integer multiples of $\scale$. The generator finally
is
\begin{equation*}
    [\gen\testf](x)=\frac{\inNvar^2}{2}\,\partial_x^2\testf(x)+
    \param\frac{\testf(x+\scale)-2\testf(x)+\testf(x-\scale)}{2}
\end{equation*}
because the integral jump term reduces itself to a finite difference
term.

\subsection{The \rqm\ laws $\relat(\param)$}\label{relqmLaw}

The family of the \rqm\ (quantum mechanics) laws on the other hand
is a particular case of the well known (centered and symmetric)
Generalized--Hyperbolic family (see for example~\cite{cufaro07} and
references quoted therein): in fact we have
$\relat(\param)=\mathfrak{GH}\left(-\frac{1}{2},1,\param\right)$, as
can be seen by direct inspection of their \PDF's and \CHF's. Remark
as a consequence that these are not simple Hyperbolic laws that
constitute the different particular sub-family
$\mathfrak{GH}\left(1,1,\param\right)$. The name follows from the
fact that -- for a suitable identification of the parameters
$\param$ and $\scale$ by means of the particle mass $m$, the
velocity of light $c$ and the Planck constant $\hbar$ -- its \pd\
generator
\begin{equation*}
   A= \lch(\partial_x)=mc^2-\sqrt{m^2c^4-c^2\hbar^2\partial^2_x}
\end{equation*}
coincides with the Hamiltonian operator of the simplest form of a free relativistic Schr\"odinger
equation~\cite{cufaro09,applebaum}, and hence its corresponding \LS\ equation exactly coincides
with this free relativistic Schr\"odinger equation (see Appendix~\ref{rqmequations}).
$\relat(\param)$ is closed under convolution, as can be seen from the form of the \CHF's, but the
laws are not stable for the same reasons as the \VG: the parameter $\param$ is not a scale
parameter. The \PDF's and \CHF's are explicitly known (see Table~\ref{families}), and all their
moments exist: the odd moments (in particular the expectation) vanish by symmetry, while the even
moments are always finite and its variance is $\param$. Since the L\évy measure is explicitly
known (see Table~\ref{families})  the L\évy dimensionless generator also takes the form
\begin{equation*}
    [\gen\testf](x)=\param\Zeroint{y}{[\testf(x+y)-\testf(x)]}{\,\frac{K_1(|y|)}{\pi |y|}}.
\end{equation*}
where $K_1$ is a modified Bessel function.

\section{The L\évy--Schr\"odinger equation}\label{LS}

To keep the notations as simple as possible also in this chapter the laws and the time coordinate
will again be supposed dimensionless. It has been shown in~\cite{cufaro09} that the evolution
equation\refeq{symmFP} of a centered, symmetric L\évy process can be formally turned into a \LS\
equation: in fact the \pd\ generator $\lch(\partial)$ of our processes is a self--adjoint operator
in $L^2$ and hence can correctly play the role of a hamiltonian. We summarize in the following the
formal steps leading to the \LS\ equation (for further details see~\cite{cufaro09}); this will
also establish the notation for the subsequent sections.

Take as  background noise a centered, symmetric, \id\ law with
$\pdf,\; \chf=e^\lch,\;\Ltriple=(0,\diff,\Lpdf)$ and a symmetric
$\Lpdf$ so that\refeq{sLKh} holds
    \begin{equation*}
        \lch(u)=-\frac{\diff^2}{2}\,u^2+\Zeroint{y}{(\cos uy
        -1)}{\,\Lpdf(y)};
    \end{equation*}
remember that since $\lch$ is real and symmetric, $\chf$ too will be
real, symmetric and non negative ($\chf\geq0$). Define then the
transition \CHF\ $\trchf(u,t)=\chf^t(u)$ and the reduced transition
\PDF\
    \begin{equation*}
        \trpdf(x,t)=\frac{1}{2\pi}\totint{\chf^t(u)e^{-iux}}{u}
    \end{equation*}
of the corresponding L\évy process, and take an initial law
$\inpdf$, $\inchf=e^{\inlch}$: the \CHF\ and the \PDF\ of the
process will  be
\begin{eqnarray}
    \prchf(u,t)&=&\trchf(u,t)\inchf(u)\nonumber\\
    \prpdf(x,t)&=&[\trpdf(t)*\inpdf](x)=\totint{\trpdf(x-y,t)\inpdf(y)}{y}\nonumber\\
    \prpdf(x,t)&=&\frac{1}{2\pi}\totint{\prchf(u,t)e^{-iux}}{u}.\label{invchf}
\end{eqnarray}
There are hence two ways to calculate $\prpdf(x,t)$: either as
$\prpdf=\trpdf*\inpdf$, or by inverting the \CHF\
$\prchf=\trchf\inchf$. As a matter of fact these two ways  give the
same result, but -- depending on the specific problem -- one can be
easier to calculate than the other. The \PDF\ $\prpdf(x,t)$ of the
previous step must also be a solution of the (dimensionless)
evolution equation
\begin{equation}\label{FPeq}
    \partial_t\prpdf(x,t)=\frac{\diff^2}{2}\partial_x^2\prpdf(x,t)+\Zeroint{y}{[\prpdf(x+y,t)-\prpdf(x,t)]}{\,\Lpdf(y)}
\end{equation}
and in principle we could find $\prpdf$ also by directly solving
this equation.

We pass then to the \LS\ propagators by means of the formal
substitution $t\rightarrow it$:
\begin{equation*}
    \propchf(u,t)=\trchf(u,it)=\chf^{it}(u)=e^{it\lch(u)},\qquad\quad
    \prop(x,t)=\trpdf(x,it)
\end{equation*}
so that $\prop$ and $\propchf$ will still verify the same
reciprocity relations\refeq{trpdf} of $\trpdf$ and $\trchf$
\begin{equation*}
    \prop(x,t)=\trpdf(x,it)=\frac{1}{2\pi}\totint{\trchf(u,it)e^{-iux}}{u}=\frac{1}{2\pi}\totint{\propchf(u,t)e^{-iux}}{u}.
\end{equation*}
Remark that if the law of the background noise is centered, symmetric and \id\ then $\lch$ is
real, symmetric and positive and hence we always have $|\propchf|=1$. This implies first that
$\propchf$ is not normalizable in $L^2$, and hence that also $\prop$ is not normalizable in $L^2$.
This is not surprising since, as it is well known, the propagators are not supposed to be
normalizable \WF's. On the other hand, as we will see later, this also entails that an initial
normalized \WF\ will stay normalized all along its evolution. We choose now an initial \LS\ \WF:
to compare the evolutions of the \WF's with that of the process \PDF's, we will start -- whwnever
we can --  with a law $\inpdf$, $\inchf=e^{\inlch}$ and with a \WF\ $\inwf$ such that
$|\inwf|^2=\inpdf$, namely
    \begin{equation*}
        \inwf(x)=\sqrt{\inpdf(x)}\,e^{iS_0(x)}
    \end{equation*}
where $S_0$ is an arbitrary, dimensionless, real function. In this
way we are also sure that $\inwf\in L^2(\Real)$, and that
$\|\inwf\|^2=1$. As a matter of fact we could also characterize our
initial state through the \WF\ \FT\ $\inhatwf(u)$ which exists
because $\inwf\in L^2(\Real)$. It is  possible to show that $\inchf$
and $\inhatwf$ must satisfy the following relation
    \begin{equation*}
        \overline{\chf}_0=\inhatwf*\hat{\overline{\wf}}_0
    \end{equation*}
which is simply the dual of $|\inwf|^2=\inpdf$. The initial \WF\ can
be simplified by choosing $\inpdf$ and $\inchf$ centered and
symmetric, with $S_0=0$. In this way we will have real $\inchf$ and
$\inwf$, with
   \begin{equation}\label{sqrtpdf}
        \inwf(x)=\sqrt{\inpdf(x)} ,
    \end{equation}
so that the following relation will always be satisfied
    \begin{equation}\label{FTconv}
        \inchf=\inhatwf*\inhatwf.
    \end{equation}
Now the \LS\ \WF's will obey the following evolution scheme
\begin{eqnarray}
    \hatwf(u,t)&=&\propchf(u,t)\inhatwf(u) \nonumber\\
    \wf(x,t)&=&[\prop(t)*\inwf](x)=\totint{\prop(x-y,t)\inwf(y)}{y}\nonumber\\
    \wf(x,t)&=&\frac{1}{\sqrt{2\pi}}\totint{\hatwf(u,t)e^{iux}}{u}\label{wfFT}
\end{eqnarray}
Here we can see the relevance of having $|\propchf|^2=1$ (namely of
having a centered, symmetric background L\évy noise, and hence a
self--adjoint generator): we have indeed that
$|\hatwf(t)|^2=|\propchf|^2|\inwf|^2=|\inwf|^2$, so that if
$\|\inwf\|^2=1$ then also $\|\hatwf(t)\|^2=1$, and as a consequence
(by Parseval and Plancherel theorems) $\|\wf(t)\|^2=1$ at every $t$.
In other words we can say that the non normalizability of the
propagator is the counterpart of the unitarity of the \LS\
evolution. Finally the \WF's $\wf(x,t)$ introduced in the previous
steps must satisfy the free \LS\ equation
    \begin{equation}\label{LSeq}
    i\partial_t\wf(x,t)=-\frac{\diff^2}{2}\partial_x^2\wf(x,t)-\Zeroint{y}{[\wf(x+y,t)-\wf(x,t)]}{\,\Lpdf(y)}.
    \end{equation}

\section{Processes and wave packets}\label{wpack}

We will give now several examples of \LS\ \WF's compared with the corresponding purely L\évy
evolutions. We classify these examples first by choosing the laws of the background noises: this
will be done by picking up the \id\ laws that allow a reasonable knowledge of both the transition
\PDF\ of the L\évy process, and the \LS\ propagator. Besides the usual Wiener case (that will be
considered just to show the way) this will indeed allow us to calculate the evolutions by means of
integrations, without being obliged to solve \pd\ equations. The equation will be used instead --
when it is possible -- as a check on the solutions found from transition \PDF's and propagators.
We will compare then the typical evolutions of the L\évy process \PDF's, and of the \WF's
solutions of a free \LS\ equation: for details, notations and formulas about both the initial laws
and \WF's, and the transition \PDF's and propagators we will make due references to
Appendix~\ref{initialstates} and to Appendix~\ref{transprop}. Remark also that in the following we
will reintroduce the dimensional parameters $\scale, \inscale$ and $\tscale$.

\subsection{Gauss}

Take a Wiener process with transition law\refeq{transnorm}: for a normal initial law\refeq{innorm}
$\norm_{\,\inscale}$ we have
\begin{equation*}
    \prchf(u,t)=\trchf(u,t)\inchf(u)=e^{-(2Dt+\inscale^2)u^2/2)}
\end{equation*}
so that the evolution is always Gaussian $\norm(2Dt+\inscale^2)$: it starts with a non degenerate
normal distribution of variance $\inscale^2$ and then widens as the usual diffusions do with
variance $2Dt+\inscale^2$. The \LS\ evolution of the \WF's on the other hand is here the usual
quantum mechanical one: take first as initial \WF\ the Gaussian\refeq{innormwf}: we then have as
wave packets
\begin{eqnarray*}
    \hatwf(u,t)&=&\propchf(u,t)\inhatwf(u)=\sqrt[4]{\frac{2\inscale^2}{\pi}}\,e^{-(\inscale^2+iDt)u^2} \\
    \wf(x,t)&=&\frac{1}{\sqrt{2\pi}}\totint{\hatwf(u,t)e^{iux}}{u}=\sqrt[4]{\frac{\inscale^2}{2\pi}}\,
                   \frac{e^{-x^2/4(\inscale^2+iDt)}}{\sqrt{\inscale^2+iDt}}
\end{eqnarray*}
It is well known that in this case $|\wf(x,t)|^2$ has a widening,
Gaussian shape all along its evolution. We neglect to display
pictures of these well known evolutions.

\subsection{Cauchy}\label{cauchyevolutions}

The Cauchy process is one of the most studies non Gaussian, L\évy processes~\cite{garbaczewski1},
first  because it is stable, and then because the calculations are relatively accessible. For
example, if the initial law is a Cauchy $\cauchy_{\,\inscale}$ with $\trchf(u,t)=e^{-\vel t|u|}$,
form\refeq{transcauchy} and\refeq{incauchy} we immediately have for the transition \CHF\
\begin{equation*}
    \prchf(u,t)=e^{-(\inscale+\vel t)|u|}
\end{equation*}
namely the process law remains a Cauchy $\cauchy_{\,\inscale+\vel
t}$ at every $t$ with a typical broadening for $t\to+\infty$
\begin{equation}\label{cauchycauchy}
    \prpdf(x,t)=\frac{1}{\pi}\,\frac{b+\vel t}{(b+\vel t)^2+x^2}.
\end{equation}
Of course this behavior (which is in common with the Gaussian Wiener
process) comes out from the fact that the Cauchy laws are stable,
and we neglect to display the corresponding figure.
\begin{figure}
\begin{center}
\includegraphics*[width=13.0cm]{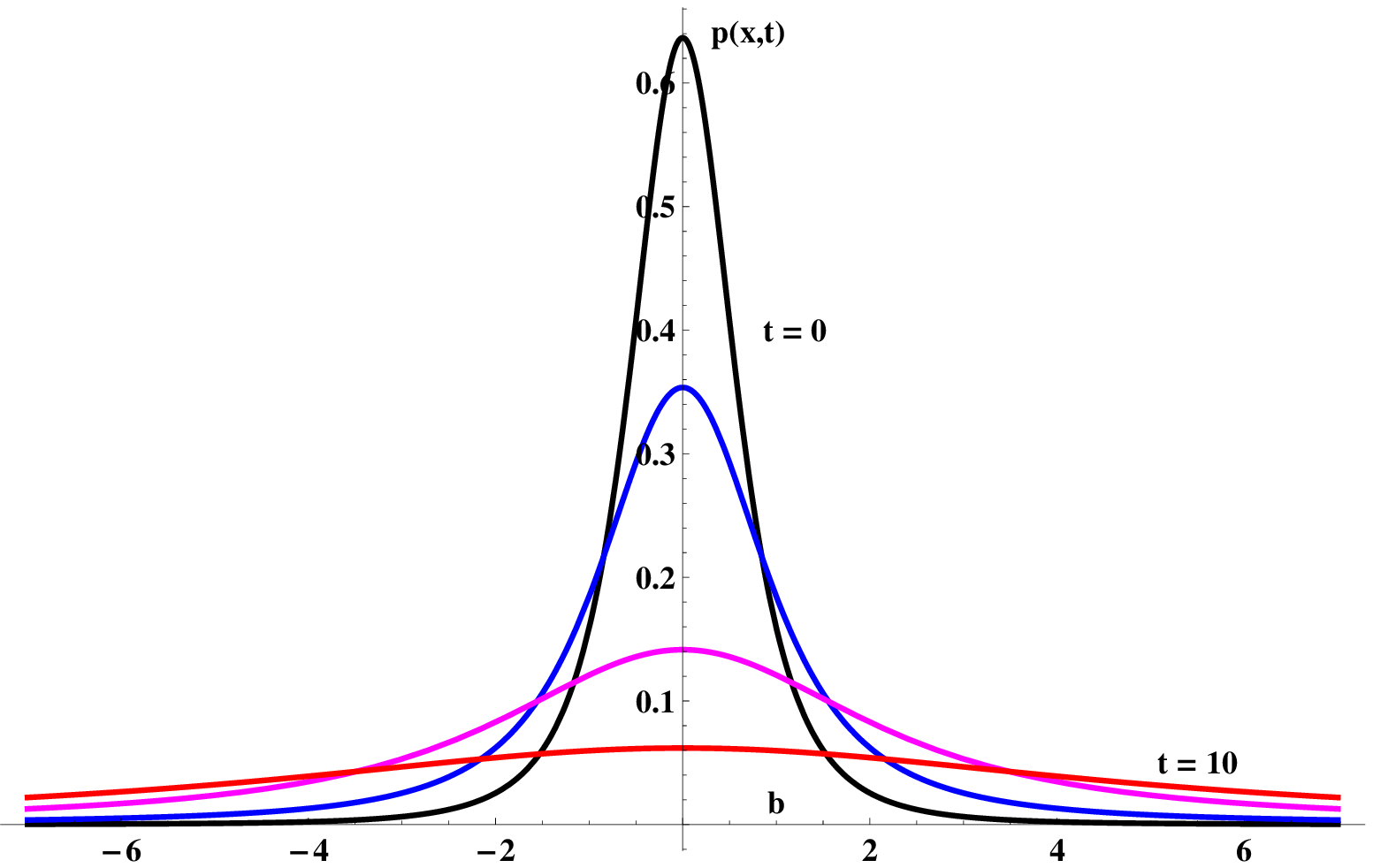}
\end{center}
\caption{The \PDF\refeq{cauchystudent} for a Cauchy process with a
Student $\stud_\inscale(3)$ initial distribution.}\label{fig04_b}
\end{figure}
Even when the initial \PDF\ is a $\stud_\inscale(3)$ with
$\inchf(u)=(1+b|u|)e^{-\inscale|u|}$ calculations are easy: now the
transition law is again $\cauchy_{\vel t}$, and the one--time
process law $\cauchy_{\vel t}*\stud_\inscale(3)$ will have as \CHF
\begin{equation*}
    \prchf(u,t)=\trchf(u,t)\inchf(u)=(1+b|u|)e^{-(\inscale+\vel t)|u|}
\end{equation*}
while the \PDF\ is recovered by \CHF\ inversion:
\begin{equation}\label{cauchystudent}
    \prpdf(x,t)=\frac{1}{2\pi}\totint{\prchf(u,t)e^{-iux}}{u}=\frac{(\inscale+\vel t)^2(2\inscale+\vel t)+vtx^2}{\pi\left[(\inscale+\vel t)^2+x^2\right]^2}.
\end{equation}
It would be easy to check that this is again a normalized,
uni--modal, bell--shaped, broadening \PDF\ (see Fig.~\ref{fig04_b}),
with neither an expectation nor a finite variance for $t>0$. We
would find in particular that the process law is the mixture
\begin{equation}\label{betamix}
  \cauchy_{\vel t}*\stud_\inscale(3) = \frac{1}{2}\,\frac{\vel t}{\inscale+\vel t}\,\widetilde{\mathfrak{B}}^{1/2}_{\inscale+\vel t}\left(\frac{3}{2},\frac{1}{2}\right)
     +\frac{1}{2}\,\frac{2\inscale+\vel t}{\inscale+\vel t}\,\widetilde{\mathfrak{B}}^{1/2}_{\inscale+\vel t}\left(\frac{1}{2},\frac{3}{2}\right)
\end{equation}
of the laws $\widetilde{\mathfrak{B}}^{1/2}_{\scale}(\alpha,\beta)$
of the square root of second kind Beta \rv's (see
Appendix~\ref{beta} for details). Remark that in particular
$\widetilde{\mathfrak{B}}^{1/2}_{\inscale+\vel
t}(1/2,3/2)=\stud_{\inscale+\vel t}(3)$. For this example we can
also show by direct calculation that the \PDF's\refeq{cauchycauchy}
and\refeq{cauchystudent} are both solutions of the \pd\ Cauchy
equation\refeq{cauchyeq}.

\begin{figure}
\begin{center}
\includegraphics*[width=13.0cm]{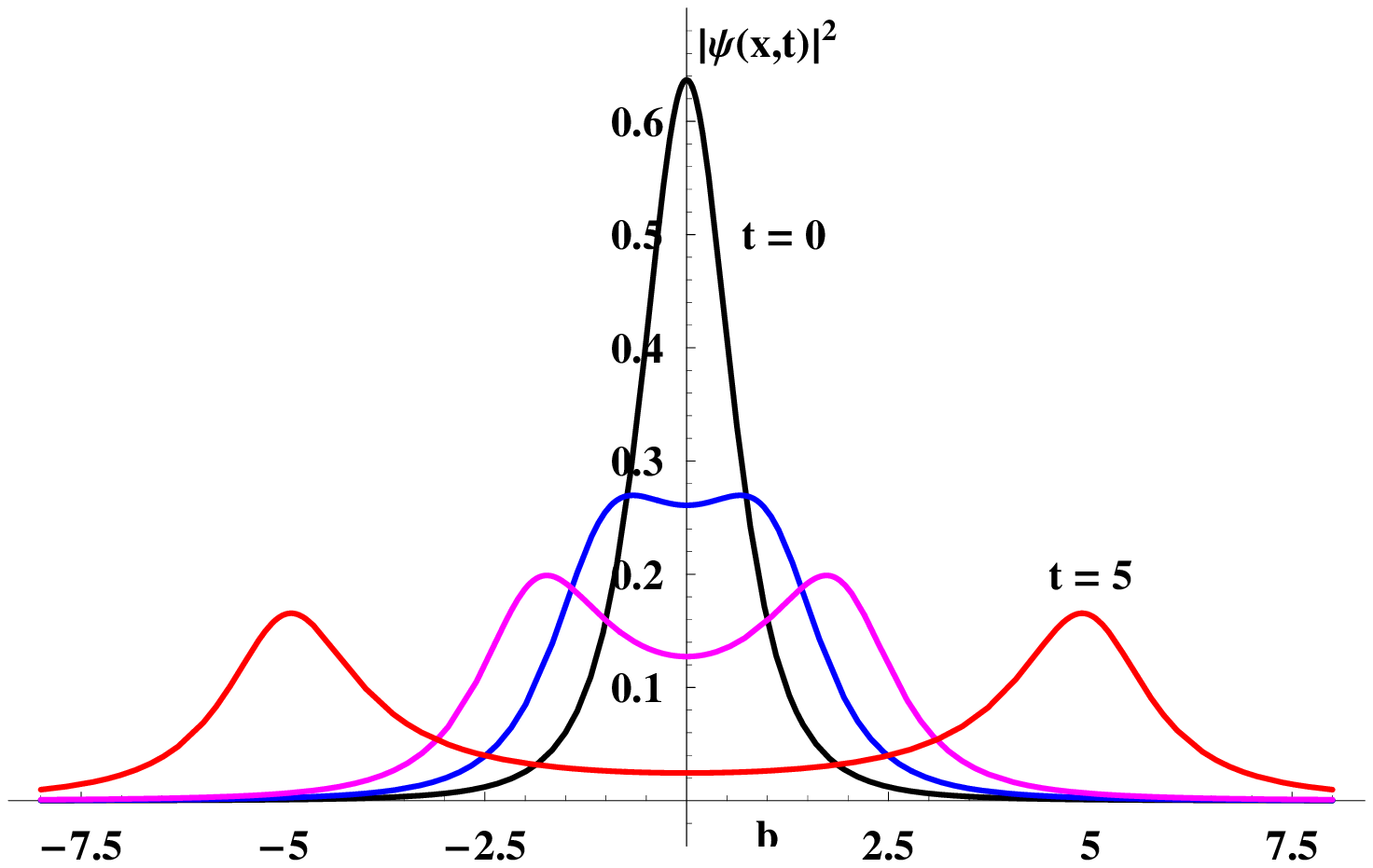}
\end{center}
\caption{The square modulus of the Cauchy--Schr\"odinger \WF\refeq{studentWF} for a Student
$\stud_\inscale(3)$ initial distribution.}\label{fig05}
\end{figure}
\begin{figure}
\begin{center}
\includegraphics*[width=13.0cm]{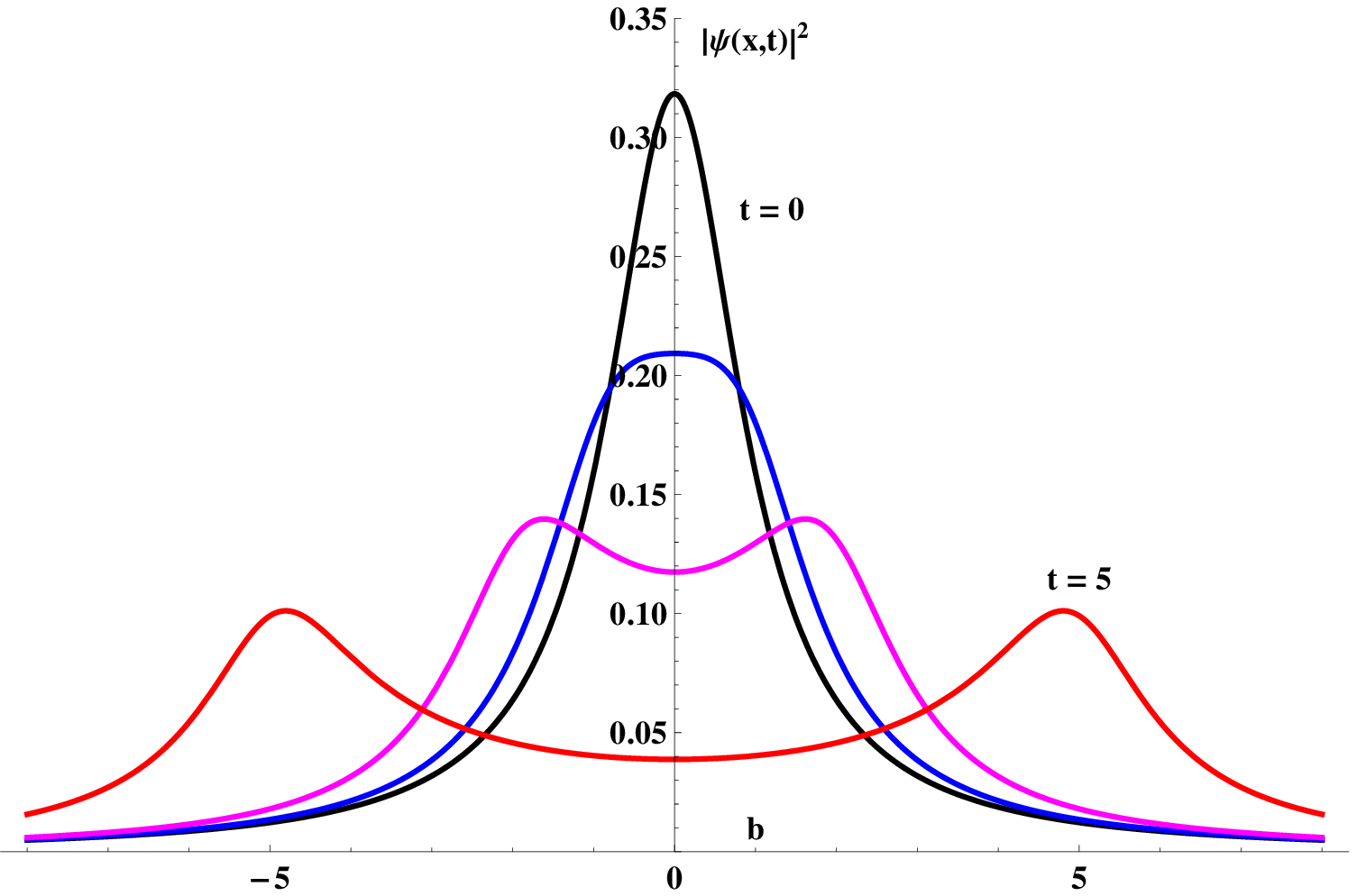}
\end{center}
\caption{The square modulus of the Cauchy--Schr\"odinger
\WF\refeq{cauchyWF} for a Cauchy $\cauchy_{\,\inscale}$ initial
distribution.}\label{fig06}
\end{figure}

The Cauchy--Schr\"odinger evolutions, on the other hand, show a more
interesting structure. The simplest case is found when we take as
$|\inwf|^2$ the Student $\stud_\inscale(3)$ case\refeq{in3studwf}:
from\refeq{propcauchy} indeed we have
\begin{equation*}
    \hatwf(u,t)=\propchf(u,t)\inhatwf(u)=\sqrt{\inscale}\,e^{-(b+i\vel t)|u|}
\end{equation*}
and hence
\begin{equation}\label{studentWF}
    \wf(x,t)=\frac{1}{\sqrt{2\pi}}\totint{\hatwf(u,t)e^{iux}}{u}
    =\sqrt{\frac{2\inscale}{\pi}}\,\frac{\inscale+i\vel t}{(\inscale+i\vel t)^2+x^2}
\end{equation}
This \WF\ (see Figure~\ref{fig05}) is correctly normalized in $L^2$
but shows a new feature: \emph{bi-modality}. In fact $|\wf|^2$ has
now two well defined maxima smoothly drifting away from the center
as $t\to+\infty$. It is also possible to show -- as an example --
that our \WF\ is a solution of the Cauchy--Schr\"odinger
equation\refeq{cauchyschreq}. For the right--hand side of this
equation we indeed have from the principal value integral
\begin{equation*}
   -\Zeroint{y}{[\wf(x+y,t)-\wf(x,t)]}{\frac{\vel}{\pi y^2}}
      =\vel\sqrt{\frac{2\inscale}{\pi}}\,\frac{(\inscale+i\vel t)^2-x^2}{\left[(\inscale+i\vel t)^2+x^2\right]^2}
\end{equation*}
which is easily seen to coincide with $i\partial_t\wf(x,t)$. As a
consequence the \WF\refeq{studentWF} correctly satisfies the \pd\
Cauchy--Schr\"odinger equation\refeq{cauchyschreq}. A similar result
is found in the case of a Cauchy $\cauchy_{\,\inscale}$ initial
\WF\refeq{incauchywf}: from the propagator\refeq{propcauchy} we have
\begin{equation*}
    \hatwf(u,t)=\frac{\sqrt{2\inscale}}{\pi}\,K_0(\inscale|u|)e^{-i\vel t|u|}
\end{equation*}
and hence by inverting the \FT:
\begin{eqnarray}
  \wf(x,t) &=&\frac{1}{\sqrt{2\pi}}\totint{\frac{\sqrt{2\inscale}}{\pi}\,K_0(\inscale|u|)e^{-i\vel
  t|u|}e^{iux}}{u}\nonumber\\
   &=& \frac{1}{\pi\sqrt{\inscale\pi}}
          \left[A\left(\frac{x+\vel t}{\inscale}\right)+\overline{A\left(\frac{x-\vel t}{\inscale}\right)}\,\,\right]\label{cauchyWF}
\end{eqnarray}
where we defined
\begin{equation*}
    A(z)=\frac{\frac{\pi}{2}-i\,\arcsinh z}{\sqrt{1+z^2}}
\end{equation*}
and we used the following two results
\begin{equation*}
    \halfint{\cos(xz)K_0(z)}{z}=\frac{\pi}{2}\,\frac{1}{\sqrt{1+x^2}}\,,
       \qquad\quad\halfint{\sin(xz)K_0(z)}{z}=\frac{\arcsinh
       x}{\sqrt{1+x^2}}\,.
\end{equation*}
The \WF\refeq{cauchyWF} is normalized in $L^2$ and shows (see
Figure~\ref{fig06}) a behavior similar to that of\refeq{studentWF}:
its \PDF\ $|\psi|^2$ starts as a Cauchy $\cauchy_{\,\inscale}$
distribution and then widens with two well defined maxima drifting
away from the center. Here too, hence, we have bi-modality: remark
the difference with the Cauchy process \PDF's
$\cauchy_{\,\inscale+\vel t}$ and $\cauchy_{\vel
t}*\stud_\inscale(3)$ which instead broaden by remaining strictly
unimodal.

\subsection{Laplace}

This bi-modality of the wave packets, or at least its breaking in two symmetric structures
drifting away from the center can also be found in other examples. Take first the \VG\ process of
Appendix~\ref{vgprocess}. At variance with the Cauchy process, this is an example of a non stable,
\sd\ process and hence has a certain interest as a non typical case. We will refer to the
Appendix~\ref{vginstate} for a discussion of possible initial states. At present we will limit our
discussion to initial states of the same \VG\ family of the background noise, and we will also
always choose coincident scale parameters $\scale=\inscale$ for the background noise and the
initial states.

\begin{figure}
\begin{center}
\includegraphics*[width=13.0cm]{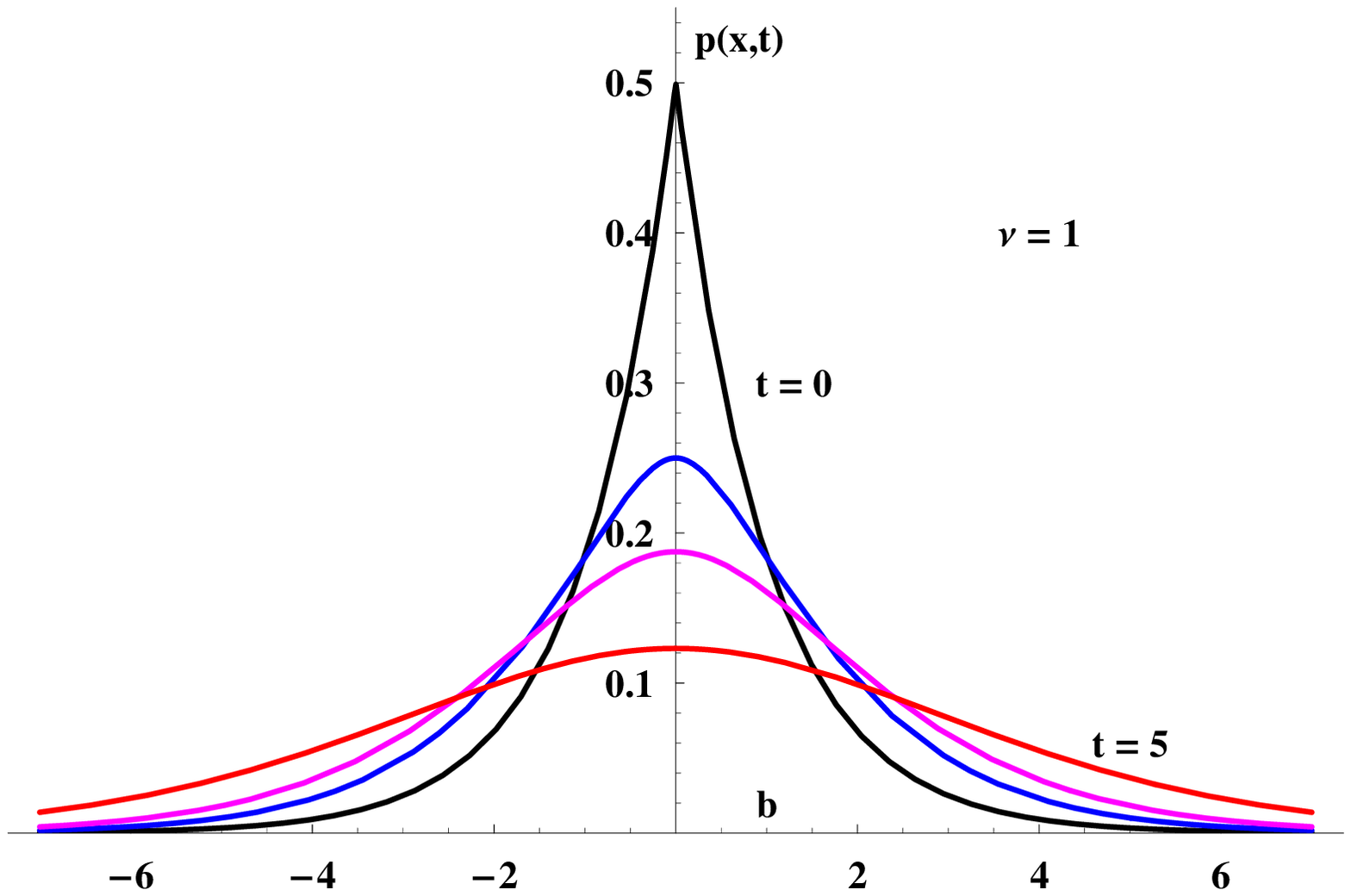}
\end{center}
\caption{The \PDF\refeq{vginvg} for a \VG\ process with Laplace $\vg_\inscale(1)=\lapl_\inscale$
initial distribution.}\label{fig08}
\end{figure}
\begin{figure}
\begin{center}
\includegraphics*[width=13.0cm]{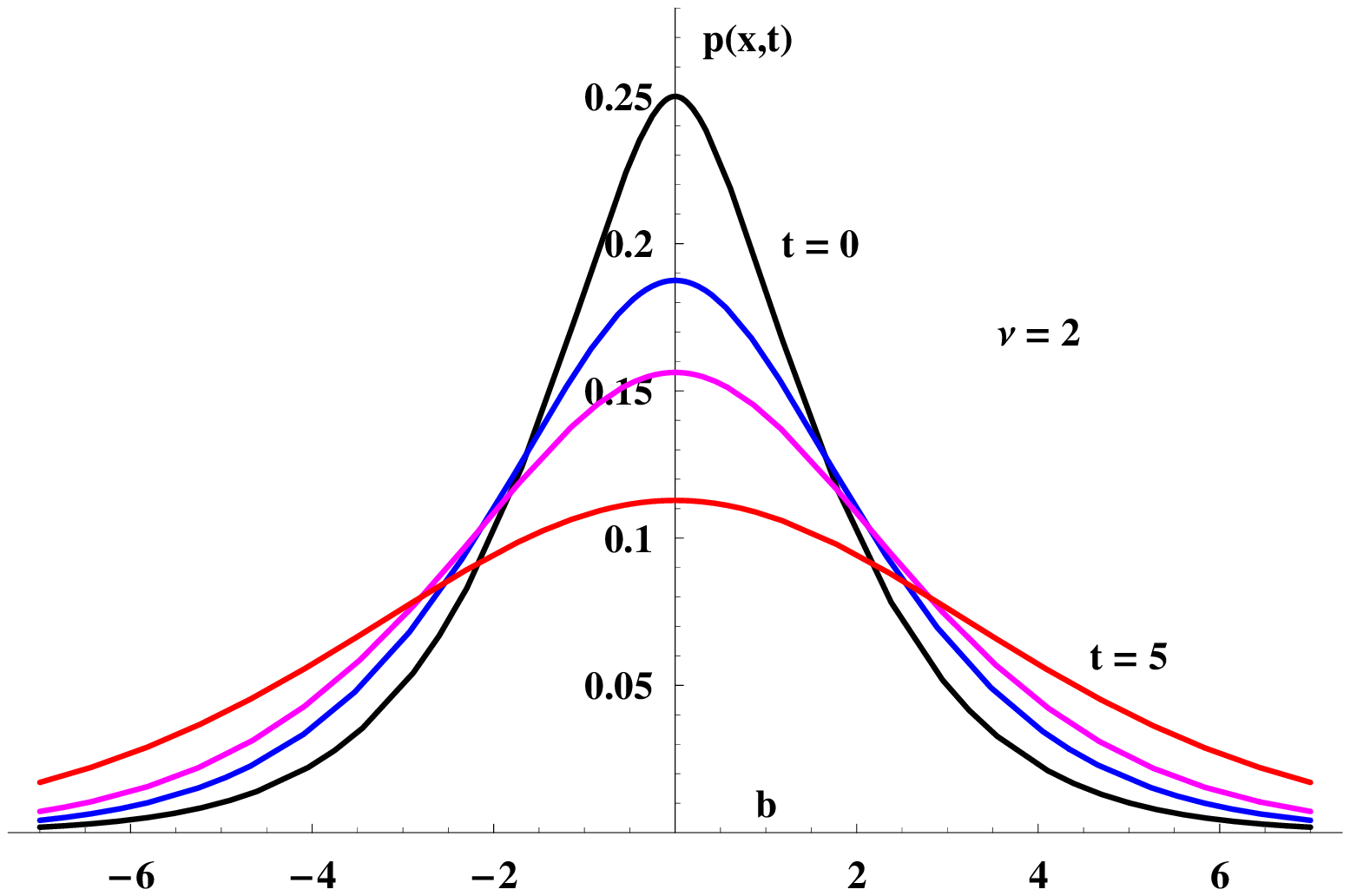}
\end{center}
\caption{The \PDF\refeq{vginvg} for a \VG\ process with \VG\ $\vg_\inscale(2)$ initial
distribution.}\label{fig09}
\end{figure}
For a \VG\ process with transition law\refeq{transvg} and initial
\PDF\refeq{invg} we immediately have
\begin{equation*}
    \prchf(u,t)=\trchf(u,t)\inchf(u)=\left(\frac{1}{1+\inscale^2u^2}\right)^{\nu+\freq t}
\end{equation*}
and hence the process law simply is $\vg_\inscale(\nu+\freq t)$ with \PDF
\begin{equation}\label{vginvg}
    \prpdf(x,t)=\frac{2}{2^\nu\Gamma(\nu)\sqrt{2\pi}\,\inscale}\left(\frac{|x|}{\inscale}\right)^{\nu+\freq t-\frac{1}{2}}\!
           K_{\nu+\freq t-\frac{1}{2}}\left(\frac{|x|}{\inscale}\right)
\end{equation}
namely always a \VG\ but with a growing parameter $\nu+\freq t$. On
the one hand this explains why it would be delusory to think of
simplifying the example by starting, for instance, with a Laplace
$\lapl_\inscale=\vg_\inscale(1)$ initial law: in fact at every time
$t>0$ the process law would in any case no longer be a Laplace law,
but a more general \VG\ with $\nu+\freq t\neq1$. On the other hand
this apparently explains why at every $t$ the \PDF\ will appear as a
broadening, uni--modal distribution as shown in the
Figures~\ref{fig08} and~\ref{fig09} respectively for $\nu=1$ and
$\nu=2$.

\begin{figure}
\begin{center}
\includegraphics*[width=13.0cm]{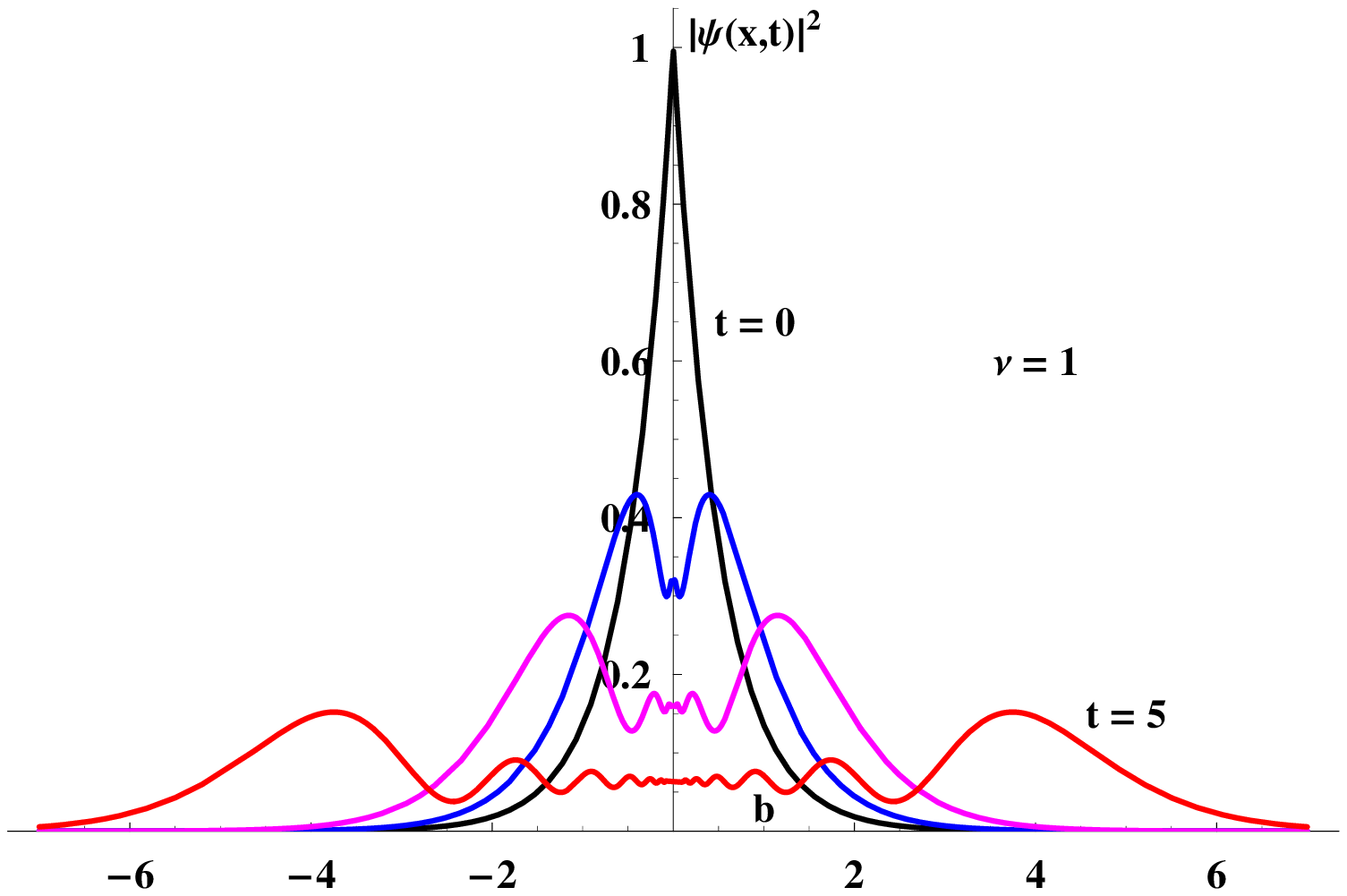}
\end{center}
\caption{The square modulus of the \VG--Schr\"odinger \WF\refeq{vgwfinvg} with a Laplace
$\vg_\inscale(1)=\lapl_\inscale$ initial \WF.}\label{fig10}
\end{figure}
\begin{figure}
\begin{center}
\includegraphics*[width=13.0cm]{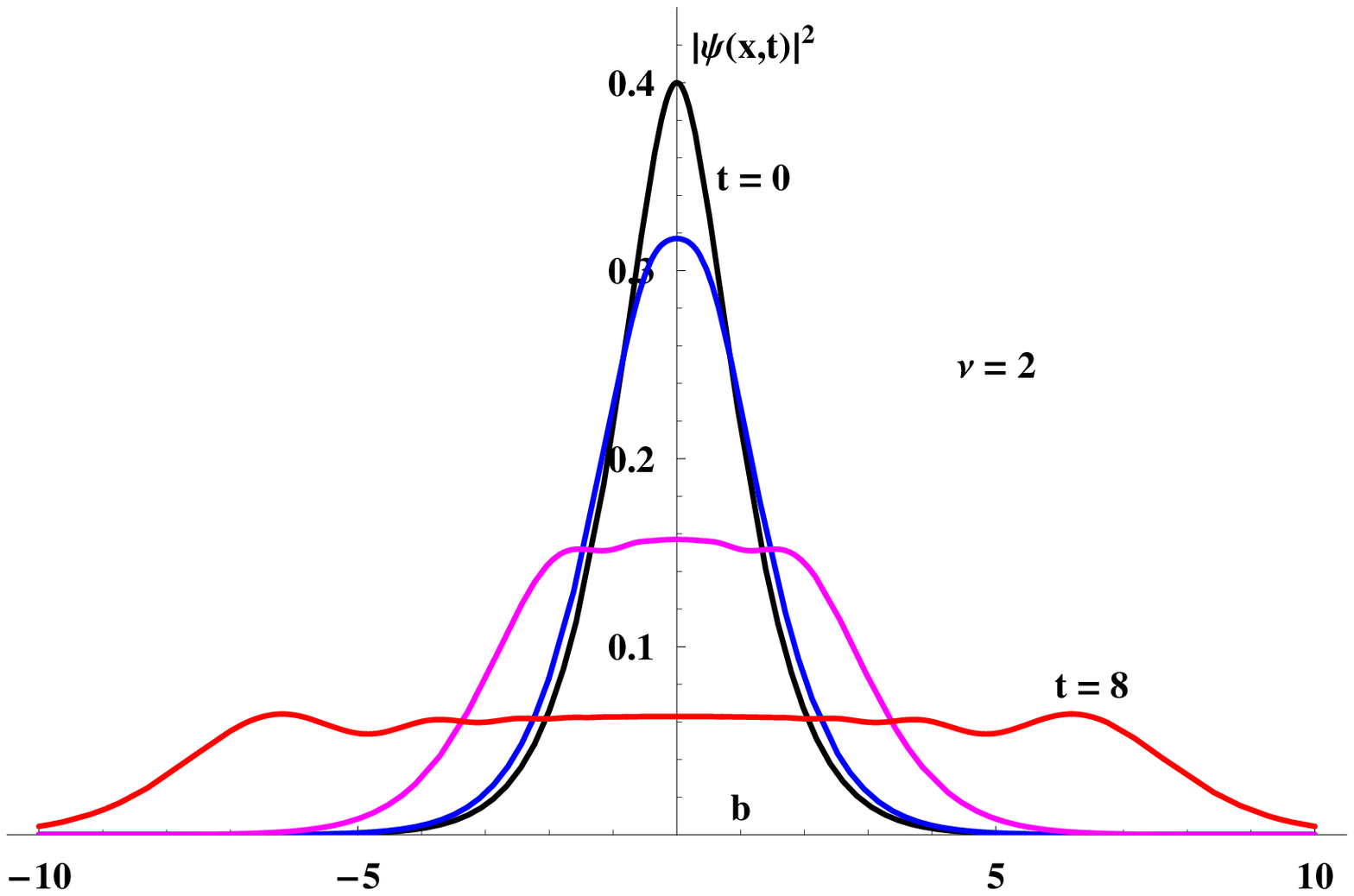}
\end{center}
\caption{The square modulus of the \VG--Schr\"odinger \WF\refeq{vgwfinvg} with a \VG\
$\vg_\inscale(2)$ initial \WF.}\label{fig11}
\end{figure}
For a \LS\ evolution, on the other hand, we have from\refeq{propvg}
and\refeq{invgwf2}
\begin{equation*}
    \hatwf(u,t)=\sqrt{\frac{\inscale}{\sqrt{\pi}}\,\frac{\Gamma(2\nu)}{\Gamma\left(2\nu-\frac{1}{2}\right)}}\left(\frac{1}{1+\inscale^2u^2}\right)^{\nu+i\freq t}
\end{equation*}
so that the inverse \FT\ will be
\begin{eqnarray}
  \wf(x,t) &=& \frac{1}{2\pi}\totint{\hatwf(u,t)e^{iux}}{u} \nonumber\\
   &=&\sqrt{\frac{\inscale}{\sqrt{\pi}}\,\frac{\Gamma(2\nu)}{\Gamma\left(2\nu-\frac{1}{2}\right)}}\,
      \frac{2}{2^{\nu+i\freq t}\Gamma(\nu+i\freq
      t)\sqrt{2\pi}} \nonumber\\
   &&\qquad \qquad\qquad\qquad\qquad  \frac{1}{\inscale}\left(\frac{|x|}{\inscale}\right)^{\nu+i\freq+1/2}
                             K_{\nu+i\freq+1/2}\left(\frac{|x|}{\inscale}\right)\label{vgwfinvg}
\end{eqnarray}
Numerical calculations and plotting then show that the
\WF\refeq{vgwfinvg} always is normalized, and that $|\wf|^2$ has two
maxima symmetrically drifting away from the center (see
Figure~\ref{fig10}). The behavior in $x=0$ is rapidly oscillating,
but with infinitesimal amplitude as we approach $x=0$: in fact the
singular behavior of the Bessel function is here competing with an
infinitesimal $|x|^\nu$ factor. The distribution shows also a slowly
decreasing, flat plateau (with micro--oscillations) in the central
region, while the diverging maxima can be rather dull as in the
Figure~\ref{fig11}.

\subsection{Poisson}

\begin{figure}
\begin{center}
\includegraphics*[width=13.0cm]{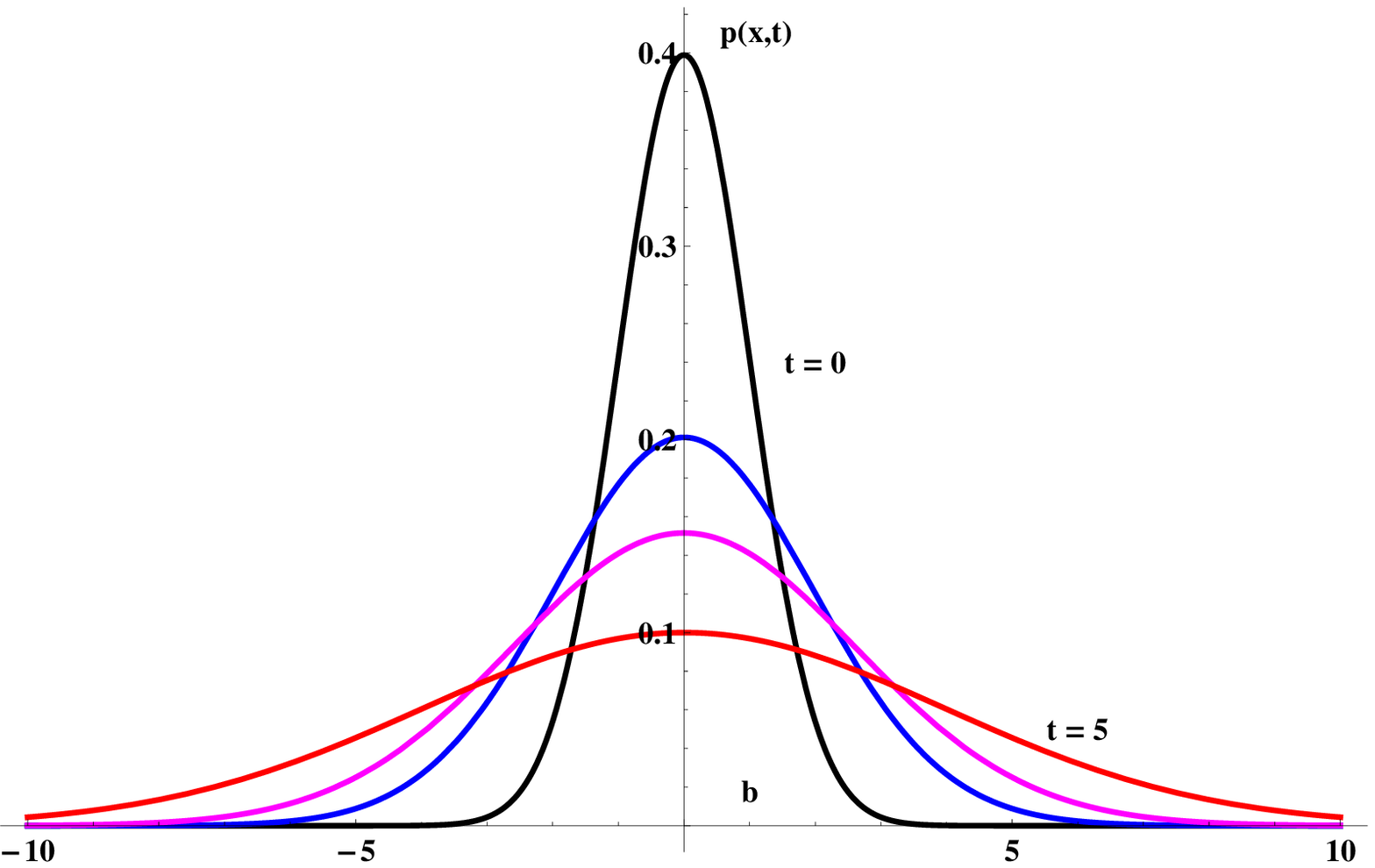}
\end{center}
\caption{The \PDF\refeq{gpd} of a Normal--Poisson process
$\norm(2Dt)*\Poiss\left(\freq t,\degen_\scale\right)$ with a
Gaussian initial law.}\label{fig12_b}
\end{figure}

The following examples  come from two \ac, but not \sd\ background noises: the compound
Wiener--Poisson processes introduced in the Appendix~\ref{comppoiss}. First take the process with
the transition law $\norm(2Dt)*\Poiss\left(\freq t,\norm_\Nvar\right)$ in\refeq{transnp}: with a
normal initial law\refeq{innorm} the marginal law of the process becomes
$\norm(2Dt+\inscale^2)*\Poiss\left(\freq t,\norm_\Nvar\right)$ namely
\begin{equation}\label{gpn}
    \prpdf(x,t)=e^{-\freq t}\sum_{k=0}^{\infty}\frac{(\freq
    t)^k}{k!}\,\frac{e^{-x^2/2(k\Nvar^2+2Dt+\inscale^2)}}{\sqrt{2\pi(k\Nvar^2+2Dt+\inscale^2)}}
\end{equation}
which apparently is a Poisson mixture of centered, normal \PDF's of
different variances, and hence has the usual bell--like, uni--modal,
diffusing shape that we will not bother to show. For the other
transition law $\norm(2Dt)*\Poiss\left(\freq t,\degen_\scale\right)$
in\refeq{transdp} with the same normal initial distribution the
marginal law instead is $\norm(2Dt+\inscale^2)*\Poiss\left(\freq
t,\degen_\scale\right)$ namely
\begin{equation}\label{gpd}
    \prpdf(x,t) = e^{-\freq t}\sum_{k=0}^{\infty}\frac{(\freq t)^k}{k!}\,\frac{1}{ 2^k}\sum_{j=0}^k\binom{k}{j}
                                \frac{e^{-[x-(k-2j)\,\scale]^2/2(2Dt+\inscale^2)}}{\sqrt{2\pi(2
                                Dt+\inscale^2)}}.
\end{equation}
In other words we always have generalized Poisson mixtures, but of
non centered normal \PDF's. Even in this case, however, the shape of
the overall \PDF\ will be that of a bell--like, uni--modal,
diffusing curve (see Figure~\ref{fig12_b}).

\begin{figure}
\begin{center}
\includegraphics*[width=13.0cm]{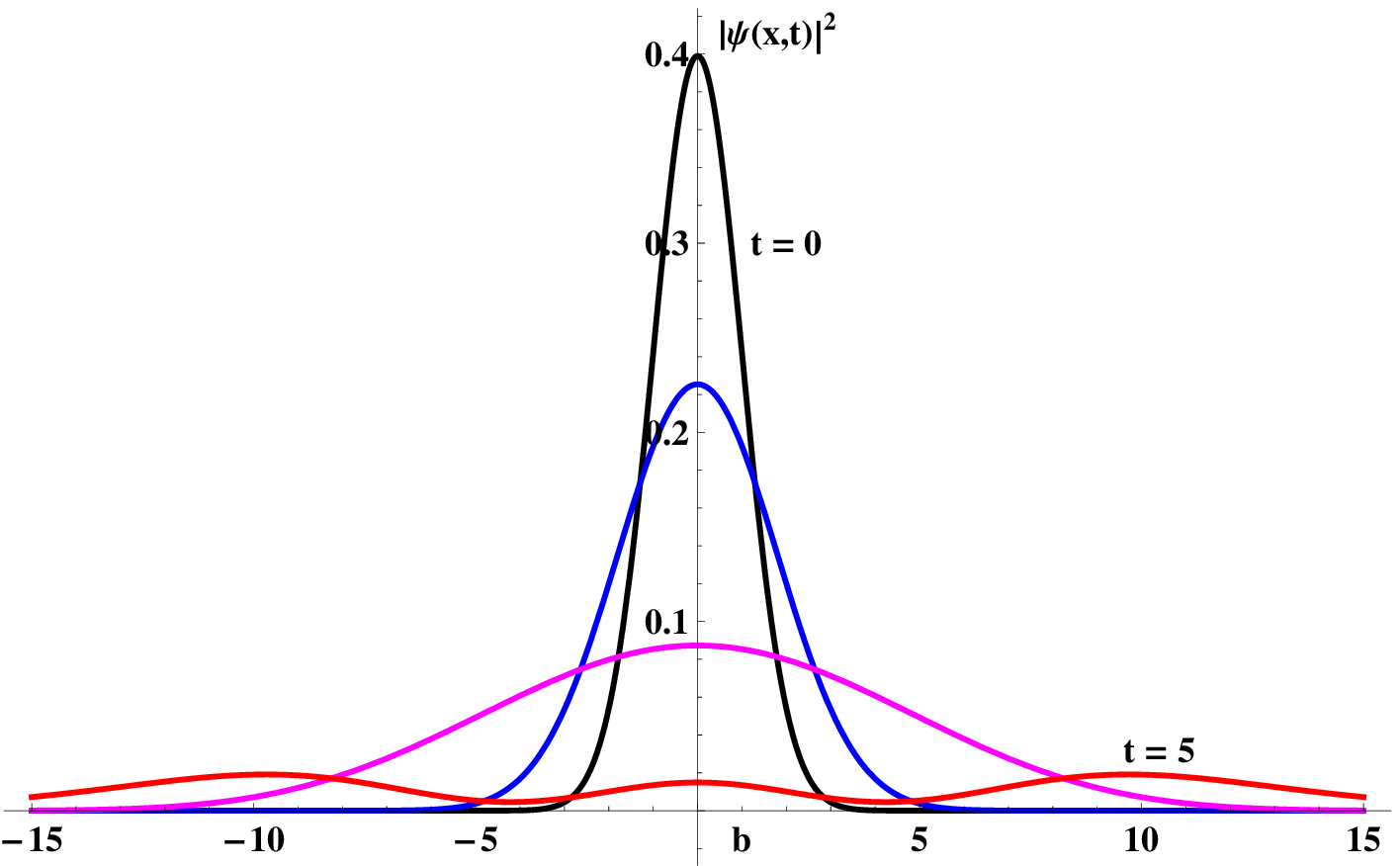}
\end{center}
\caption{The square modulus of the Normal--Poisson Schr\"odinger
\WF\refeq{lsgpd} with a Gaussian initial \WF.}\label{fig12}
\end{figure}

For the \LS\ equation on the other hand consider first the
propagator $\norm\left(2iDt\right)*\Poiss\left(i\freq
t,\norm_\Nvar\right)$ in\refeq{propnp} applied to an initial
Gaussian \WF\refeq{innormwf}; we then have
\begin{equation*}
   \hatwf(u,t)=e^{i\freq t\left(e^{-\Nvar^2u^2/2}-1\right)}\sqrt[4]{\frac{2\inscale^2}{\pi}}\,\,e^{-(\inscale^2+iDt)u^2}
\end{equation*}
and, by inverting the \FT\ and taking into account the properties of the Gaussian integrals, the
\WF\ will be
\begin{equation}\label{lsgpn}
    \wf(x,t)=e^{i\freq t}\sum_{k=0}^{\infty}\frac{(i\freq t)^k}{k!}\sqrt[4]{8\pi\inscale^2}\,\frac{e^{-x^2/2(k\Nvar^2+2\inscale^2+2iDt)}}{\sqrt{2\pi(k\Nvar^2+2\inscale^2+2iDt)}}
\end{equation}
namely a time--dependent, complex, Poisson superposition of centered
Gaussian \WF's. The same is true for the second example with
propagator $\norm\left(2iDt\right)*\Poiss\left(i\freq
t,\degen_\scale\right)$ in\refeq{propdp} with an initial Gaussian
\WF\refeq{innormwf}: the \WF\ \FT\ in fact now is
\begin{equation*}
   \hatwf(u,t)=e^{i\freq t\left(\cos\scale u-1\right)}\sqrt[4]{\frac{2\inscale^2}{\pi}}\,\,e^{-(\inscale^2+iDt)u^2}
\end{equation*}
so that the \WF\ itself will be
\begin{equation}\label{lsgpd}
    \wf(x,t)=e^{i\freq t}\sum_{k=0}^{\infty}\frac{(i\freq
    t)^k}{k!}\frac{\sqrt[4]{8\pi\inscale^2}}{2^k}\sum_{j=0}^k\binom{k}{j}
    \frac{e^{-[x-(k-2j)\scale]^2/4(\inscale^2+iDt)}}{\sqrt{4\pi(\inscale^2+iDt)}}.
\end{equation}
In conclusion, while the plots of $\prpdf(x,t)$ in\refeq{gpn}
and\refeq{gpd} simply display the too familiar story of a diffusing
bell--shaped curve, and the same would be true for $|\wf(x,t)|^2$
in\refeq{lsgpn}, for $|\wf(x,t)|^2$  in\refeq{lsgpd} we instead have
again a separation of the wave packet in two symmetrical
sub--packets drifting away from the center (see Figure~\ref{fig12}).

\subsection{\Rqm}
\begin{figure}
\begin{center}
\includegraphics*[width=13.0cm]{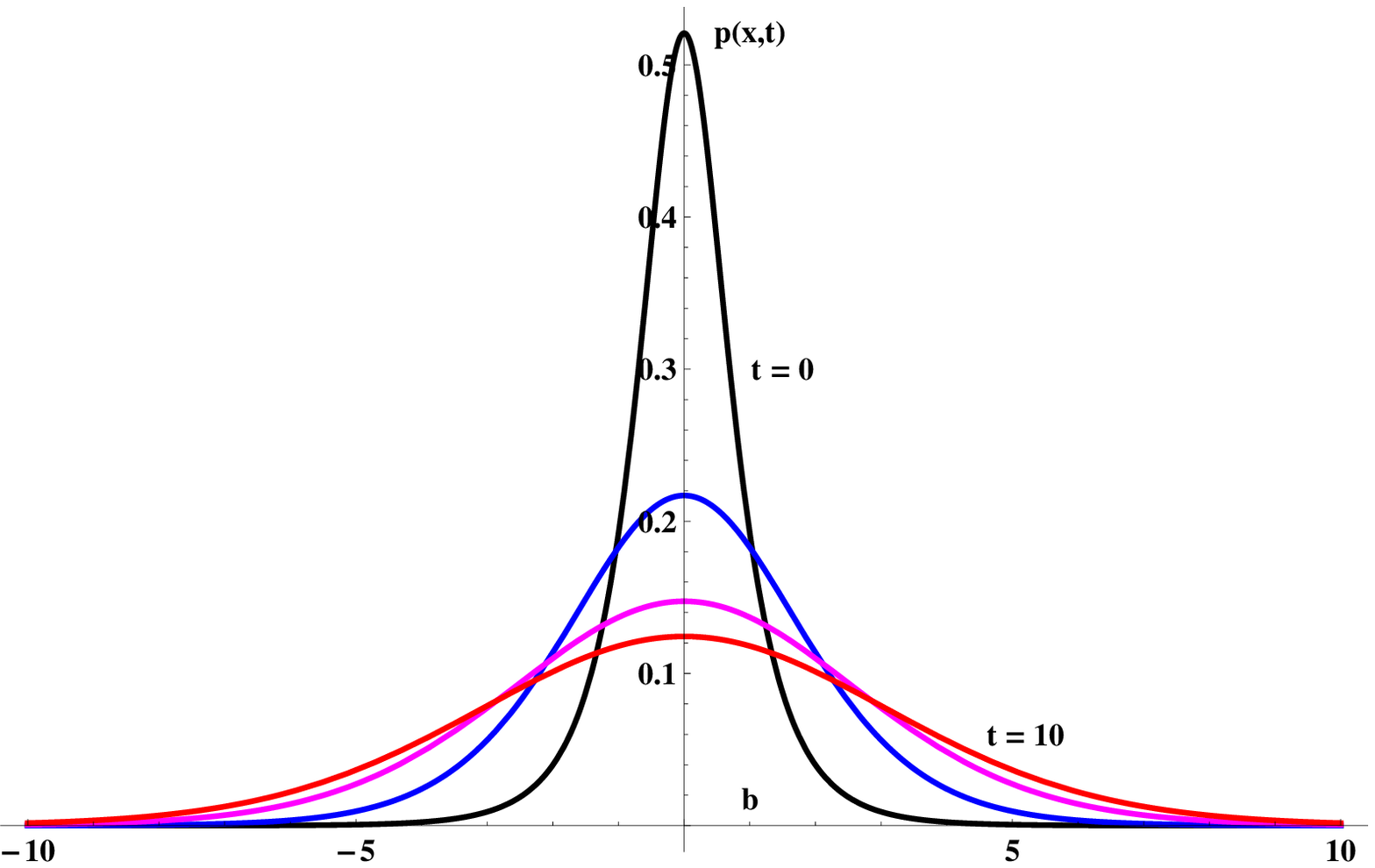}
\end{center}
\caption{The \PDF\refeq{rqmpdf} of a L\évy process with a \Rqm\
background noise and an initial law of the same
family.}\label{fig13}
\end{figure}

\begin{figure}
\begin{center}
\includegraphics*[width=13.0cm]{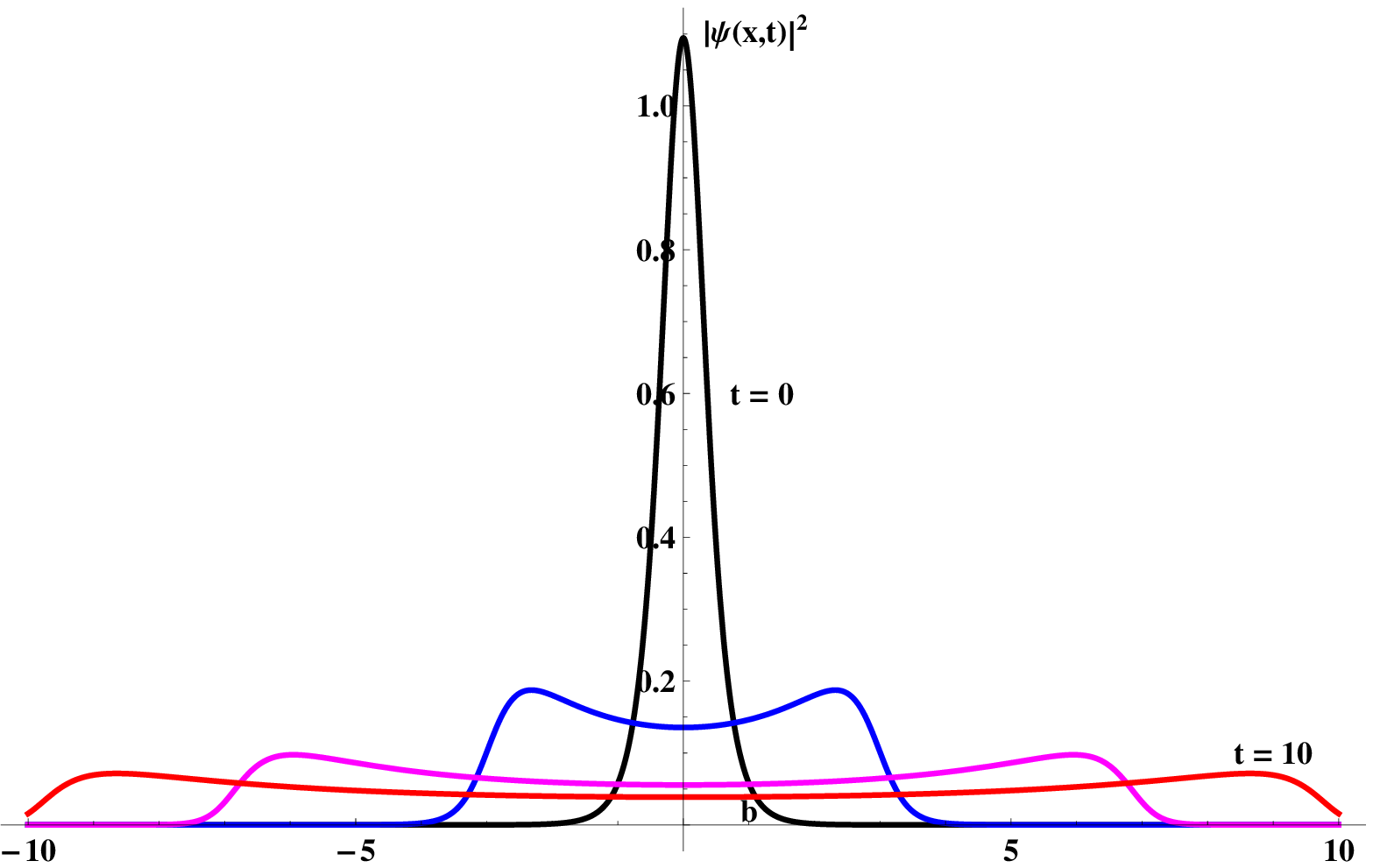}
\end{center}
\caption{The square modulus of the \Rqm\ \WF\refeq{rqmwp} with an
initial \WF\ of the same family.}\label{fig14}
\end{figure}

In a way similar to that of the \VG, for a \Rqm\ L\évy process with
transition law\refeq{rqmtr} and initial distribution\refeq{rqmin},
but with $\scale=\inscale$, we immediately have
\begin{eqnarray}
    \prchf(u,t)&=&\trchf(u,t)\inchf(u)=e^{(\nu+\freq
    t)\left(1-\sqrt{1+\scale^2u^2}\right)}\label{rqmchf}\\
    \prpdf(x,t)&=&\frac{(\nu+\freq t)e^{\nu+\freq t}}{\pi\scale}\,\,\frac{K_1\left(\sqrt{(\nu+\freq t)^2+x^2/\scale^2}\right)}{\sqrt{(\nu+\freq
    t)^2+x^2/\scale^2}}\label{rqmpdf}
\end{eqnarray}
and hence the process law simply is $\relat(\nu+\freq t)$, namely it
will stay always in the same \Rqm\ family but with a time dependent
parameter. The \PDF\ $\prpdf(x,t)$ is shown in the
Figure~\ref{fig13} and has the usual bell--like, uni--modal,
diffusing form. For the corresponding \LS\ evolution on the other
hand we have from\refeq{rqminwf} and\refeq{rqmprop} that the
normalized \WF's are
\begin{eqnarray}
    \hatwf(u,t)&=&\propchf(u,t)\inhatwf(u)=\sqrt{\frac{\scale}{2e^{2\nu}K_1(2\nu)}}\,\,e^{(\nu+i\freq
      t)\left(1-\sqrt{1+\scale^2u^2}\right)}\label{rqmwp1}\\
    \wf(x,t)&=& \frac{(\nu+i\freq t)e^{i\freq t}}{\sqrt{\scale\pi K_1(2\nu)}}\,\,
      \frac{K_1\left(\sqrt{(\nu+i\freq t)^2+x^2/\scale^2}\right)}{\sqrt{(\nu+i\freq t)^2+x^2/\scale^2}}\label{rqmwp}
\end{eqnarray}
We show in the Figure~\ref{fig14} how this $|\wf(x,y)|^2$ behaves,
and in particular, at variance with the previous L\évy
\PDF\refeq{rqmpdf}, we find here again that the the \WF\ shows two
symmetric maxima drifting away from the center of the distribution:
the bi-modality that we have already pointed out in all our other
\LS\ examples.

\section{Conclusions}\label{conclusions}

We presented in the previous sections several examples of free wave
packets that are solutions of the \LS\ equation without
potentials\refeq{LSeq}. We started by generalizing the relation
between Brownian motion and Schr\"odinger equation, and by
associating the kinetic energy of a physical system to the generator
of a symmetric L\'evy process, namely to a pseudo-differential
operator whose symbol is the \LCH\ $\lch$ of an \id\ law. This
amounts to suppose, then, that the \LS\ equation is based on an
underlying L\'evy process that can have both Gaussian (continuous)
and non Gaussian (jumping) components. The use of all the \id, even
non stable, processes on the other hand is important and physically
meaningful because there are significant cases that are in the
domain of our \LS\ picture, without being in that of the stable
(fractional) Schr\"odinger equation. In particular, as discussed
in~\cite{cufaro09,garbaczewski1}, the simplest form of a
relativistic, free Schr\"odinger equation can be associated with a
particular type of \sd, non stable process acting as background
noise. Moreover in many instances of the L\'evy--Schr\"odinger
equation the new energy--momentum relations can be seen as
corrections to the classical relations for small values of certain
parameters~\cite{cufaro09}. It must also be remembered that -- at
variance with the stable, fractional case -- our model is not tied
to the use of processes with infinite variance: the variances can be
chosen to be finite even in a purely non Gaussian model -- as in the
case of the relativistic, free Schr\"odinger equation -- and can
then be used as a legitimate measure of the dispersion. Finally let
us recall that a typical non stable, Student L\'evy noise seems to
be suitable for applications in the models of halo formation in
intense beam of charged particles in
accelerators~\cite{applications,cufaro07,vivoli}.

It was then important to explore the general behavior of the diffusing \LS\ \WF's: we
systematically approached this problem by defining in Section~\ref{LS} a procedure allowing us to
explore several combinations of initial \WF's (Appendix~\ref{initialstates}) and background L\évy
noises (Appendix~\ref{transprop}), and by comparing L\évy processes and free \LS\ wave packets. We
have then remarked that virtually in all our examples of Section~\ref{wpack} we witnessed a
similar qualitative behavior: first of all the \LS\ wave packets diffuse, in the sense that they
broaden in a very regular way. As it is known the variance of a L\évy process -- when it exists --
grows linearly with the time, exactly as in the usual diffusions. Of course stable, non Gaussian
noises are excluded, since for them there is no variance, and we have instead an anomalous sub-
and super-diffusive behavior. The corresponding \LS\ wave packets show a similar qualitative
behavior also if it is not always easy to calculate their variances.

A second, more surprising feature however is represented by the bi-modality of the \LS\ \WF's. In
fact we found that in virtually all our examples the wave packet splits in two sub--packets
symmetrically and smoothly drifting away from the center: a behavior which is present neither in
the free L\évy processes, nor in the (Gaussian) free Schr\"odinger \WF's. It is interesting to
remark, then, that the unique instance with a similar bi--modal behavior has been found
earlier~\cite{chechkin} deals with \emph{confined} L\évy flights. In our opinion the bi--modality
found in our examples could then be connected to the combined effect of Nelson dynamics, and L\évy
jumps in the background noise, and it would be interesting to explore if this behavior shows up
again in form of rings and shells respectively for the two- and three-dimensional \LS\ equation.
This bi-modality, on the other hand, is in sheer contrast with the uni--modality of both the L\évy
processes and the (Gaussian) Schr\"odinger \WF's.

It would be important now to explicitly give in full detail the formal association between \LS\
\WF's and the underlying L\évy processes, namely a true generalized stochastic mechanics. In
particular we would show that to every \WF\ solution of the \LS\ equation we can associate a well
defined L\évy process: the techniques of the stochastic calculus applied to L\'evy processes are
today in full development~\cite{sato,applebaum,protter}, and at our knowledge there is no
apparent, fundamental impediment along this road. Finally it would be relevant to explore this
L\évy--Nelson stochastic mechanics by adding suitable potentials to our \LS\ equation, and by
studying the corresponding possible stationary and coherent states: all that too will be the
subject of future papers.

\appendix

\section{Types of laws}\label{types}

As stated in the Section~\ref{notation} we deal in this paper with centered laws of \rv's $X$.
Even when the expectation does not exist we can always speak of centering around the
\textit{median}. On the other hand to eliminate the centering it will be enough to take $X+b$ with
$b\in\Real$ instead of $X$, then to substitute $x-b$ to $x$ in the $\pdf$, and to add a factor
$e^{ibu}$ to the \CHF\ $\chf$. For our purposes it will also be expedient to introduce a
dimensional \emph{scale parameter} $\scale>0$ to take into account the physical dimensions of our
\rv's: to fix the ideas in this paper $a$ will be supposed to be a \textit{length}. Take first a
\rv\ $X$ with law $\law$, \PDF\ $\pdf$ and \CHF\ $\chf$, and suppose that $X$ is a dimensionless
quantity; then the variables argument of $\pdf$ and $\chf$, will be dimensionless. On the other
hand $X_\scale=\scale X$ will be a length and will follow a law $\law_\scale$ with
\begin{equation*}
    \pdf_\scale(x)\,dx=\pdf\left(\frac{x}{\scale }\right)\frac{dx}{\scale }\,,\qquad\quad
    \chf_\scale(u)=\chf(\scale u).
\end{equation*}
Here $x$ and $u$ are now dimensional variables ($x$ is a length,
while $u$ is the reciprocal of a length), so that $x/\scale$ and
$\scale u$ will be dimensionless. Remark that within this notation
we numerically have $\law=\law_1$, so that for instance
$\pdf_1(x)=\pdf(x)$. This could be slightly misleading since the
argument of $\pdf_1$ is a length, while that of $\pdf$ is supposed
to be dimensionless. To avoid any possible misunderstanding we will
then reserve the symbols $\law$, $\pdf$ and $\chf$ for the
dimensionless laws, while $\law_1$, $\pdf_1$ and $\chf_1$ will be
associated to the dimensional ones. For example if $X$ follows the
\textit{standard}, dimensionless normal law $\norm$ with
\begin{equation*}
    \pdf(x)=\frac{e^{-x^2/2}}{\sqrt{2\pi}}\,,\qquad\quad\chf(u)=e^{-u^2/2}
\end{equation*}
the dimensional \rv's $X_\scale=\scale X$ will follow the laws $\norm_\scale=\norm(\scale^2)$ with
\begin{equation*}
    \pdf_\scale(x)=\frac{e^{-x^2/2\scale ^2}}{\scale \sqrt{2\pi}}\,,\qquad\quad\chf_\scale(u)=e^{-\scale ^2u^2/2}.
\end{equation*}
Then $\pdf$ and $\pdf_1$ will be coincident, but the dimensional meaning of their respective
variables will be different. Remark finally that in general we will choose dimensionless laws that
are not necessarily standard laws: of course (when the variances exist) we will have
$\var{X_\scale}=\scale ^2\,\var{X}$, but $\var{X}$ is not always supposed to be equal to 1.

We could now think to $\law_\scale$ as the parametric family of the
rescaled \rv's $\scale X$: these parametric families spanned just by
one scale parameter $\scale$ are here entire \textit{types of
laws}\footnote{A \emph{type of laws} (see~\cite{loeve} Section 14)
is a family of laws that only differ among themselves by a centering
and a rescaling: in other words, if $\varphi(u)$ is the \CHF\ of a
law, all the laws of the same type have \CHF's
$e^{ibu}\varphi(\scale u)$ with a centering parameter
$b\in\mathbb{R}$, and a scaling parameter $a>0$ (we exclude here the
sign inversions). In terms of \rv's this means that the laws of $X$
and $aX+b$ (for $a>0$, and $b\in\mathbb{R}$) always are of the same
type, and on the other hand that $X$ and $Y$ belong to the same type
if and only if it is possible to find $a>0$, and $b\in\mathbb{R}$
such that $Y$ and $aX+b$ have the same law, namely $Y\eqd aX+b$.}:
in fact, since here we only deal with centered laws (see
Section~\ref{notation}), no centering parameter $b$ is required, and
our types are spanned by means of the scale parameter $\scale$ only.
In this paper we will also consider other parametric families of
laws with some dimensionless parameter $\param$, which will not in
general be coincident with the scale parameter $\scale$. We could
then have two--parameters families $\law_\scale(\param)$, and in
general we are interested in finding which sets are closed under
convolution (namely under addition of the corresponding independent
\rv's). When a type of laws is closed under convolution (as in the
normal case of the previous example) its laws are said to be
\emph{stable}: the convolution would produce another law of the same
type, namely a law with only a different \emph{scale} parameter (in
our notation: same $\param$, but different $\scale$). If instead the
convolution produces a law of the same family, but not of the same
type (different $\param$), then the family is closed under
convolution, but its laws are not stable: this is the case, among
others, of the \VG\ laws $\vg_\scale(\param)$. Finally, when the
result of a convolution is a law not belonging at all to the family,
then $\law_\scale(\param)$ is not even closed under convolution, as
for the Student $\stud_\scale(\param)$ family.

\section{Symmetric and \ac, compound Poisson laws}\label{poisson}

Among the \id, non \sd\ laws the Poisson case stands as the most
important example, but the simple Poisson law is neither symmetric,
nor \ac. We will then generalize it in order to avoid these
shortcomings. A Poisson law $\Poiss(\param)$ is a non symmetric, non
\sd, non \ac, \id\ law without Gaussian component ($\diff=0$). The
probability is concentrated on the integer numbers with the usual
Poisson distribution so that formally
\begin{equation*}
    \pdf(x)=\sum_{k=0}^{\infty}e^{-\param}\frac{\param^k}{k!}\,\delta_{k}(x)\,,\qquad\quad
    \chf(u)=e^{\param(e^{iu}-1)}\,,\qquad\quad
    \Lpdf(x)=\param\,\delta_1(x)\,.
\end{equation*}
Both expectation and variance have value $\param$. Since
$\Poiss(\param)$ is neither centered, nor symmetric the generator of
the corresponding L\évy process will not be self--adjoint. It is
well known, moreover, that the sample paths of the corresponding
simple Poisson process are ascending staircase trajectories, with
randomly located steps of unit height, $\param$ representing the
average number of jumps per unit time interval. As a consequence
these processes are not \ac. To move ahead we must then first
symmetrize the Poisson law, and then make it \ac.

Take a symmetric (we do not require it to be \ac\ or \id) law $\complaw$ with \CHF\
$\compchf(u)=e^{\complch(u)}$ and build the corresponding compound Poisson law
$\Poiss(\param,\complaw)$ with \CHF
\begin{equation*}
    \chf(u)=e^{-\param}\sum_{k=0}^{\infty}\frac{\param^k}{k!}\,\compchf^k(u)=e^{\param[\compchf(u)-1]}\,,
    \qquad\quad\lch(u)=\param[\compchf(u)-1]
\end{equation*}
thus generalizing the simple Poisson case where $\compchf(u)=e^{iu}$. When $\complaw$ is also \ac\
with \PDF\ $\comppdf(x)$ the law of $\Poiss(\param,\complaw)$ is
\begin{equation}\label{cpoiss}
    \pdf(x)=e^{-\param}\sum_{k=0}^{\infty}\frac{\param^k}{k!}\,\comppdf^{*k}(x)\,,\qquad\quad
    \comppdf^{*k}=\left\{
                    \begin{array}{ll}
                      \overbrace{\comppdf*\ldots*\comppdf}^{\tiny\hbox{$k$ times}}, & \hbox{$k=1, 2 \ldots$} \\
                      \delta_0, & \hbox{$k=0$}
                    \end{array}
                  \right.
\end{equation}
but we can immediately see that this is still not \ac\ even if $\complaw$ has a density: in fact
for $k=0$ we always have a degenerate law $\delta_0$. The compound Poisson law
$\Poiss(\param,\complaw)$ has neither a drift ($\drift=0$ because of the required symmetry) nor a
Gaussian part ($\diff=0$), and its L\évy \PDF\ (that we will suppose for simplicity to show no
singularities at $x=0$) is $\Lpdf(x)=\param\comppdf(x)$: namely we have
$\Ltriple=(0,0,\param\comppdf)$. The laws of the increments of the corresponding compound Poisson
process $\Poiss(\freq t,\complaw)$ with $\freq=\param/\tscale$ are then the time dependent
mixtures
\begin{equation}\label{transpoisson}
    \prpdf(x,t)=e^{-\freq t}\sum_{k=0}^{\infty}\frac{(\freq t)^k}{k!}\,\comppdf^{*k}(x)
\end{equation}
while its self--adjoint generator (no singularities are present at $x=0$) is
\begin{equation*}
    [\gen\testf](x)=\param\totint{[\testf(x+y)-\testf(x)]\comppdf(y)}{y}\,.
\end{equation*}
Its sample trajectories are now up and down staircase functions, with steps at Poisson random
times, and random jump heights distributed according to the symmetric law $\complaw$. Since
however for $k=0$ the law is degenerate in $x=0$, these sample trajectories stick at $x=0$ for a
finite time (with probability 1), and the marginal distribution of the process is not \ac. In
other L\évy processes instead (as the Wiener process for example) the trajectory starts at $x=0$,
but its random path immediately leaves this position.

To give a first example of these symmetric (but not \ac) compound
Poisson laws take $\complaw=\norm_\Nvar$ so that
$\comppdf^{*k}\sim\norm(k\Nvar^2)$ for $k=0,1,\ldots;$ we then have
for $\Poiss(\param,\norm_\Nvar)$
\begin{equation*}
    \pdf(x)=e^{-\param}\sum_{k=0}^{\infty}\frac{\param^k}{k!}\,\frac{e^{-x^2/2k\Nvar^2}}{\sqrt{2\pi k\Nvar^2}}\,,\qquad
    \lch(u)=\param\left(e^{-\Nvar^2u^2/2}-1\right),\qquad
    \Lpdf(x)=\param\,\frac{e^{-x^2/2\Nvar^2}}{\sqrt{2\pi\Nvar^2}}\,.
\end{equation*}
The transition \PDF's of the corresponding compound Poisson process are then the time dependent
mixtures of $\norm(k\Nvar^2)$ laws
\begin{equation*}
    \prpdf(x,t)=e^{-\freq t}\sum_{k=0}^{\infty}\frac{(\freq t)^k}{k!}\,\frac{e^{-x^2/2k\Nvar^2}}{\sqrt{2\pi
    k\Nvar^2}}\,,
\end{equation*}
and the generator takes the form
\begin{equation*}
    [\gen\testf](x)=\param\totint{[\testf(x+y)-\testf(x)]\frac{e^{-y^2/2\Nvar^2}}{\sqrt{2\pi\Nvar^2}}}{y}\,.
\end{equation*}
As another example suppose instead that $\complaw=\degen_\scale$ is a Bernoulli symmetric law,
doubly degenerate around the positions $\pm \scale$, namely
\begin{equation*}
    \comppdf(x)=\frac{1}{2\scale}\,\left[\delta_1\left(\frac{x}{\scale}\right)+\delta_{-1}\left(\frac{x}{\scale}\right)\right],\qquad\compchf(u)=\cos \scale u,
\end{equation*}
and remark that now
\begin{equation*}
    \comppdf^{*k}(x)=\frac{1}{2^k\scale}\sum_{j=0}^k\binom{k}{j}\delta_{k-2j}\left(\frac{x}{\scale}\right).
\end{equation*}
As a consequence we will have for $\Poiss(\param,\degen_\scale)$:
\begin{eqnarray*}
  \pdf(x) &=& e^{-\param}\sum_{k=0}^{\infty}\frac{\param^k}{k!}\,\frac{1}{ 2^k\scale}\sum_{j=0}^k\binom{k}{j}\delta_{k-2j}\left(\frac{x}{\scale}\right) \\
  \lch(u) &=& \param(\cos \scale u-1) \\
  \Lpdf(x) &=&\frac{\param}{2\scale}\,\left[\delta_1\left(\frac{x}{\scale}\right)+\delta_{-1}\left(\frac{x}{\scale}\right)\right]
\end{eqnarray*}
We then easily have for the transition law of the process
\begin{equation*}
    \prpdf(x,t)=e^{-\freq t}\sum_{k=0}^{\infty}\frac{(\freq
    t)^k}{k!}\frac{1}{2^k\scale}\sum_{j=0}^k\binom{k}{j}\delta_{k-2j}\left(\frac{x}{\scale}\right)
\end{equation*}
while the generator is
\begin{equation*}
    [\gen\testf](x)=\frac{\param}{2}\left[\testf(x+1)-2\testf(x)+\testf(x-1)\right]
\end{equation*}

We will then further generalize our compound Poisson distributions
in order to get \ac\ laws and processes. Take a compound Poisson law
$\Poiss(\param,\complaw)$, and another independent, symmetric, \ac,
\id\ law $\complaw_0$ with \PDF\ $\backpdf(x)$, \CHF\
$\backchf(u)=e^{\backlch(u)}$ and L\évy triplet
$\Ltriple_0=(0,\diff_0,\Lpdf_0)$. Consider then the law
$\complaw_0*\Poiss(\param,\complaw)$ obtained by addition
(convolution) so that
\begin{equation*}
    \chf(u)=\backchf(u) e^{\param(\compchf(u)-1)}\,,\qquad\quad
    \lch(u)=\backlch(u)+\param(\compchf(u)-1)
\end{equation*}
while the \PDF\ is
\begin{equation*}
    \pdf(x)=e^{-\param}\sum_{k=0}^{\infty}\frac{\param^k}{k!}\,(\backpdf*\comppdf^{*k})(x)\,,\qquad
    \backpdf*\comppdf^{*k}=\left\{
                    \begin{array}{ll}
                      \backpdf*\overbrace{\comppdf*\ldots*\comppdf}^{\tiny\hbox{$k$ times}}, & \hbox{$k=1, 2 \ldots$} \\
                      \backpdf, & \hbox{$k=0$}
                    \end{array}
                  \right.
\end{equation*}
This is now a mixture of \ac\ laws. The law $\complaw_0*\Poiss(\param,\complaw)$ will also be
symmetric if both $\comppdf$ and $\backpdf$ are symmetric, and it will have a Gaussian component
if $\diff_0\neq0$. As a consequence we will have $\drift=0$ from the symmetry,
$\Lpdf(x)=\param\comppdf(x)+\Lpdf_0(x)$, and finally
$\Ltriple=\left(0\,,\,\diff_0\,,\,\param\comppdf+\Lpdf_0\right)$. The laws of the increments of
the corresponding L\évy process will then be $\chf(t)=\backchf^{t/\tscale}e^{\param
t(\compchf-1)/\tscale}$, namely
\begin{equation*}
    \lch(u,t)=\frac{t}{\tscale}\,\backlch(u)+\freq t\,[\compchf(u)-1]
\end{equation*}
so that the process will be  the superposition of two independent processes: an
$\complaw_0$--L\évy process plus a $\Poiss(\freq t,\,\complaw)$ compound Poisson process. Its
trajectories will then be the paths of the $\complaw_0$--L\évy process, interspersed with Poisson
random jumps with size law $\complaw$. If then $\backpdf(x,t)$ is the \PDF\ of
$\backchf^{\,t/\tscale}(u)$, the $t$--increment \PDF's of our process will be
\begin{equation*}
    \prpdf(x,t)=e^{-\freq t}\sum_{k=0}^{\infty}\frac{(\freq
    t)^k}{k!}\,\left[\backpdf(t)*\comppdf^{*k}\right](x)
\end{equation*}
and the self--adjoint process generator
\begin{equation*}
    [\gen\testf](x)=\frac{\diff_0^2}{2}\,\partial_x^2\testf(x)+\Zeroint{y}{[\testf(x+y)-\testf(x)]}{[\param\comppdf(y)+\Lpdf_0(y)]}\,.
\end{equation*}
Possible examples of these $\complaw_0$ are both the Gaussian and
the non Gaussian stable laws (in particular the Cauchy process), and
several self--decomposable laws as the Student or the \VG. The
relevant particular case of a Gaussian $\complaw_0$ is discussed in
the Section~\ref{poisslaw}.

\section{Initial states}\label{initialstates}

We define here a list of possible initial \PDF's and \WF's. To simplify our calculations we will
choose the initial \PDF's to be centered and symmetric, and whenever convenient we will take pairs
$\inpdf,\;\inwf$ satisfying the relation $\inpdf=|\inwf|^2$. Remark that, while $\inpdf$ is a
normalized (in $L^1$) \PDF\ and $\inchf$ is a (non normalized, and possibly non normalizable)
\CHF, $\inwf$ and $\inhatwf$ must be both normalized (in $L^2$) \WF's so that we must always pay
attention to the constants which are in front of them. Here moreover -- to put in evidence the
meaning of the involved quantities -- our laws and time coordinates will be dimensional: the space
$\scale,\inscale$ and time $\tscale$ scaling parameters will be explicitly taken into account.

\subsection{Normal $\norm_{\,\inscale}$}
Initial laws and \WF's with $\inpdf=|\inwf|^2$ are in this case
\begin{eqnarray}
    &&\inpdf(x)=\frac{e^{-x^2/2\inscale^2}}{\sqrt{2\pi \inscale^2}},\quad\qquad\inchf(u)=e^{-\inscale^2u^2/2}\label{innorm}\\
    &&\inwf(x)=\frac{e^{-x^2/4\inscale^2}}{\sqrt[4]{2\pi \inscale^2}},
    \quad\qquad\inhatwf(u)=\sqrt[4]{\frac{2\inscale^2}{\pi}}\,e^{-\inscale^2u^2}\label{innormwf}
\end{eqnarray}
Remark that, while $\inwf$ is just the square root of $\inpdf$, $\inhatwf$ is the \FT\ of $\inwf$
and its relation to $\inchf$ is given by the equation\refeq{FTconv}. The two \WF's, moreover, are
both normalized in $L^2$.

\subsection{Cauchy $\cauchy_{\,\inscale}=\stud_{\,\inscale}(1)$}
Initial laws and \WF's in this case are
\begin{eqnarray}
    &&\inpdf(x)=\frac{1}{\inscale\pi}\,\frac{\inscale^2}{\inscale^2+x^2},\qquad\qquad\inchf(u)=e^{-\inscale|u|}\label{incauchy}\\
    &&\inwf(x)=\frac{1}{\sqrt{\inscale\pi}}\sqrt{\frac{\inscale^2}{\inscale^2+x^2}},
    \,\qquad\inhatwf(u)=\frac{\sqrt{2\inscale}}{\pi}\,K_0(\inscale|u|)\label{incauchywf}
\end{eqnarray}
where $K_0$ is the modified Bessel function of order 0; it is easy to show indeed that (see for
example~\cite{abramowitz} 9.6.21)
\begin{equation*}
  \inhatwf(u) = \frac{1}{\sqrt{\inscale\pi}}\,\frac{1}{\sqrt{2\pi}}\totint{\sqrt{\frac{\inscale^2}{\inscale^2+x^2}}e^{-iux}}{x}
  =\frac{\sqrt{2\inscale}}{\pi}\,K_0(\inscale|u|).
\end{equation*}
The normalization $\|\inhatwf\|^2=1$, and the relation $\inchf=\inhatwf*\inhatwf$, are then
\begin{eqnarray*}
 \Zeroint{u}{K_0^2(\inscale|u|)}{} &=& \frac{\pi^2}{2\inscale} \\
 \int_{v\neq0,u}K_0(\inscale|u-v|)K_0(\inscale|v|)\,dv &=& \frac{\pi^2}{2\inscale}\,e^{-\inscale|u|}
\end{eqnarray*}
The first can be reduced to
\begin{equation*}
    \halfint{K_0^2(u)}{u}=\frac{\pi^2}{4}
\end{equation*}
which can be verified by direct calculation. On the other hand the convolution, that can be
reduced to the dimensionless relation
\begin{equation*}
    \int_{v\neq0,u}K_0(|u-v|)K_0(|v|)\,dv =
    \frac{\pi^2}{2}\,e^{-|u|},
\end{equation*}
does not seem to be an otherwise known result.

\subsection{3--Student $\stud_{\,\inscale}(3)$}
Initial laws and \WF's in this case are
\begin{eqnarray}
    &&\inpdf(x)=\frac{2}{\inscale\pi}\,\left(\frac{\inscale^2}{\inscale^2+x^2}\right)^{\!2},\quad\qquad\,\inchf(u)=e^{-\inscale|u|}(1+\inscale|u|)\label{in3stud}\\
    &&\inwf(x)=\sqrt{\frac{2}{\inscale\pi}}\,\frac{\inscale^2}{\inscale^2+x^2}\,,
    \qquad\qquad\inhatwf(u)=\sqrt{\inscale}\,e^{-\inscale|u|}\label{in3studwf}
\end{eqnarray}
It is  very easy to show that $\inhatwf$ is the right \FT\ of $\inwf$
\begin{equation*}
\inhatwf(u) =
\sqrt{\frac{2}{\inscale\pi}}\,\frac{1}{\sqrt{2\pi}}\totint{\frac{\inscale^2}{\inscale^2+x^2}\,e^{-iux}}{x}
= \sqrt{\inscale}\,e^{-\inscale|u|}
\end{equation*}
while here again an elementary calculation shows also that $\inchf=\inhatwf*\inhatwf$.

\subsection{\VG\ $\vg_{\,\inscale}(\nu)$}\label{vginstate}

In the general \VG\ case, to make calculations possible, we will not always choose pairs of
initial \PDF's and \WF's satisfying $\inwf=\sqrt{\inpdf}$. A possible example then is
\begin{eqnarray}
    &&\inpdf(x)=\frac{2}{2^\nu\Gamma(\nu)\sqrt{2\pi}\,\inscale}\left(\frac{|x|}{\inscale}\right)^{\nu-\frac{1}{2}}\!
           K_{\nu-\frac{1}{2}}\left(\frac{|x|}{\inscale}\right),\quad\,\inchf(u)=\left(\frac{1}{1+\inscale^2u^2}\right)^{\nu}\label{invg}\\
    &&\inwf(x)=\sqrt{\frac{2\Gamma\left(\nu+\frac{1}{2}\right)}{\inscale\pi\Gamma(\nu)\Gamma\left(2\nu-\frac{1}{2}\right)}}
    \left(\frac{|x|}{\inscale}\right)^{\nu-\frac{1}{2}}\!
           K_{\nu-\frac{1}{2}}\left(\frac{|x|}{\inscale}\right),\label{invgwf1}\\
    &&\inhatwf(u)=\sqrt{\frac{\inscale\,\Gamma(2\nu)}{\sqrt{\pi}\Gamma\left(2\nu-\frac{1}{2}\right)}}\left(\frac{1}{1+\inscale^2u^2}\right)^{\nu}\label{invgwf2}
\end{eqnarray}
where the functions are chosen  in order to have an evolution easy to calculate. The \WF's $\inwf$
and $\inhatwf$, in any case, are both normalized in $L^2$ (as can be seen by direct calculation)
and are apparently in the \FT\ relation. As a consequence here the L\évy and the \LS\ evolutions
will possibly start with different \PDF's. In fact the usual relation $\inpdf=|\inwf|^2$ could be
easily restored just in the particular case of $\nu=1$, namely for an initial Laplace law
$\lapl_{\,\inscale}=\vg_{\,\inscale}(1)$:
\begin{eqnarray}
    &&\inpdf(x)=\frac{e^{-|x|/b}}{2b},\quad\qquad\;\inchf(u)=\frac{1}{1+b^2u^2}\label{inlapl}\\
    &&\inwf(x)=\frac{e^{-|x|/2b}}{\sqrt{2b}},\quad\qquad\inhatwf(u)=\sqrt{\frac{b}{\pi}}\,\frac{2}{1+4b^2u^2}\label{inlaplwf}
\end{eqnarray}
Here it is elementary to check indeed that $\inwf=\sqrt{\inpdf}$, that $\inhatwf$ is the \FT\ of
$\inwf$, and finally that $\inchf=\inhatwf*\inhatwf$. This particular case, however, is not really
easier than the general case of the \VG\ process. In fact, as we will see soon, the parameter
affected by the time evolution is exactly $\nu$, so that it is of no help to start with $\nu=1$ if
it immediately becomes $\nu\neq1$.

\subsection{\Rqm\ $\relat_\inscale(\nu)$}

Again to make calculations easy we will choose as initial \CHF\ and \WF\ \FT\ respectively
\begin{eqnarray}
    &&\inpdf(x)=\frac{\nu e^\nu K_1\left(\sqrt{\nu^2+x^2/\inscale^2}\right)}{\inscale\pi\sqrt{\nu^2+x^2/\inscale^2}},
           \quad\qquad\inchf(u)=e^{\nu(1-\sqrt{1+\inscale^2u^2})}\label{rqmin}\\
    &&\inwf(x)=\frac{\nu K_1\left(\sqrt{\nu^2+x^2/\inscale^2}\right)}{\sqrt{\pi\inscale K_1(2\nu)(\nu^2+x^2/\inscale^2)}},
           \qquad\inhatwf(u)=\sqrt{\frac{\inscale}{2K_1(2\nu)}}\,e^{-\nu\sqrt{1+\inscale^2u^2}}\label{rqminwf}
\end{eqnarray}
which are in a relation similar to that of\refeq{invg}-(\ref{invgwf2}). It is easy to recognize
that the \WF's are correctly normalized in $L^2$.


\section{Transition laws and propagators}\label{transprop}

We will list here a few examples of background L\évy noises by paying attention to pick up
processes with a known transition \PDF\ associated to the evolution equation\refeq{FPeq} and a
known propagator associated to the free \LS\ equation\refeq{LSeq}.

\subsection{Normal $\norm(2Dt)$}

Here the background noise is a Wiener process: take a $\norm_{\,\scale}$ law with L\évy triplet
$\Ltriple=(0,\scale,0)$
\begin{equation*}
    \pdf(x)=\frac{e^{-x^2/2\scale^2}}{\sqrt{2\pi \scale^2}},\quad\qquad\chf(u)=e^{-\scale^2u^2/2}
\end{equation*}
The transition law of the corresponding L\évy process is then $\norm(2Dt)$ with
$D=\scale^2/2\tscale$, namely
\begin{equation}\label{transnorm}
    \trpdf(x,t)=\frac{e^{-x^2/4Dt}}{\sqrt{4\pi Dt}},\quad\qquad\trchf(u,t)=e^{-Dtu^2}
\end{equation}
and the \PDF\ evolution equation\refeq{FPeq} is the usual Fokker--Planck equation
\begin{equation*}
    \partial_t\prpdf(x,t)=D\,\partial^2_x\prpdf(x,t)
\end{equation*}
The corresponding \LS\ propagator $\norm(2iDt)$ is again formally normal albeit with an imaginary
variance:
\begin{equation}\label{propnorm}
    \prop(x,t)=\frac{e^{-x^2/4iDt}}{\sqrt{4\pi
    iDt}},\quad\qquad\propchf(u,t)=e^{-iDtu^2}
\end{equation}
and hence the \LS\ equation\refeq{LSeq} is the usual free Schr\"odinger equation
\begin{equation*}
    i\partial_t\wf(x,t)=-D\,\partial^2_x\wf(x,t).
\end{equation*}

\subsection{Cauchy $\cauchy_{\,\vel t}$}

From the Cauchy law $\cauchy_{\,a}$, a typical stable, non Gaussian law with L\évy triplet
$\Ltriple=(0,0,\scale/\pi x^2)$ and with
\begin{equation*}
    \pdf(x)=\frac{1}{\scale\pi}\,\frac{\scale^2}{\scale^2+x^2},\quad\qquad\chf(u)=e^{-\scale|u|}
\end{equation*}
we get the transition law $\cauchy_{\,\vel t}$ of the Cauchy process with $\vel=\scale/\tscale$:
\begin{equation}\label{transcauchy}
    \trpdf(x,t)=\frac{1}{\pi \vel t}\,\frac{\vel^2t^2}{\vel^2t^2+x^2},\quad\qquad\trchf(u,t)=e^{-\vel t|u|}
\end{equation}
and the corresponding process equation\refeq{FPeq}
\begin{equation}\label{cauchyeq}
    \partial_t\prpdf(x,t)=\Zeroint{y}{[\prpdf(x+y,t)-\prpdf(x,t)]}{\frac{\vel}{\pi
    y^2}}.
\end{equation}
On the other hand the \LS\ propagator $\cauchy_{\,i\vel t}$ is
\begin{equation}\label{propcauchy}
    \prop(x,t)=\frac{1}{i\pi}\,\frac{\vel t}{\vel^2t^2-x^2},\quad\qquad\propchf(u,t)=e^{-i\vel t|u|}.
\end{equation}
and the \LS\ equation\refeq{LSeq}
\begin{equation}\label{cauchyschreq}
    i\partial_t\wf(x,t)=-\Zeroint{y}{[\wf(x+y,t)-\wf(x,t)]}{\frac{\vel}{\pi y^2}}
\end{equation}
Remark that, at variance with the transition \PDF\refeq{transcauchy}, the Cauchy--Schr\"odinger
propagator\refeq{propcauchy} has two simple poles in $x=\pm \vel t$ drifting away from the center
$x=0$ with velocity $\vel$.

\subsection{\VG\ $\vg_\scale(\freq t)$}\label{vgprocess}

Take a \sd, non stable \VG\ law $\vg_\scale(\param)$ with symmetric L\évy triplet
$\Ltriple=(0,0,\param e^{-|x|/\scale}/|x|)$ and with
\begin{equation*}
    \pdf(x)=\frac{2}{2^\param\Gamma(\param)\sqrt{2\pi}\,\scale}\left(\frac{|x|}{\scale}\right)^{\param-\frac{1}{2}}\!
           K_{\param-\frac{1}{2}}\left(\frac{|x|}{\scale}\right),\quad\qquad\chf(u)=\left(\frac{1}{1+\scale^2u^2}\right)^{\param}
\end{equation*}
The transition law will then be $\vg_\scale(\freq t)$ with $\freq=\param/\tscale$:
\begin{equation}\label{transvg}
    \trpdf(x,t)=\frac{2}{2^{\freq t}\Gamma(\freq t)\sqrt{2\pi}\,\scale}\left(\frac{|x|}{\scale}\right)^{\freq t-\frac{1}{2}}\!
           K_{\freq t-\frac{1}{2}}\left(\frac{|x|}{\scale}\right),\quad\trchf(u,t)=\left(\frac{1}{1+\scale^2u^2}\right)^{\freq t}
\end{equation}
and the corresponding process equation\refeq{FPeq}
\begin{equation*}
    \partial_t\prpdf(x,t)=\freq\Zeroint{y}{[\prpdf(x+y,t)-\prpdf(x,t)]}{\frac{e^{-|y|/\scale}}{|y|}}.
\end{equation*}
so that the evolution will only affect the parameter $\param$, while $\scale$ will always be the
same. Then for the \LS\ propagator $\vg_\scale(i\freq t)$ we have
\begin{equation}\label{propvg}
    \prop(x,t)=\frac{2}{2^{i\freq t}\Gamma(i\freq t)\sqrt{2\pi}\,\scale}\left(\frac{|x|}{\scale}\right)^{i\freq t-\frac{1}{2}}\!
           K_{i\freq t-\frac{1}{2}}\left(\frac{|x|}{\scale}\right),\quad\propchf(u,t)=\left(\frac{1}{1+\scale^2u^2}\right)^{i\freq t}
\end{equation}
while the \LS\ equation\refeq{LSeq} becomes
\begin{equation*}
    i\partial_t\wf(x,t)=-\freq\Zeroint{y}{[\wf(x+y,t)-\wf(x,t)]}{\frac{e^{-|y|/\scale}}{|y|}}
\end{equation*}

\subsection{Wiener--Poisson
$\norm(2Dt)*\Poiss\left(\freq t,\complaw\right)$}\label{comppoiss}

We will consider here two examples of \id, non \sd\ background noise: for notations and details
see Section~\ref{poisslaw} and Appendix~\ref{poisson}. Take first the law
$\norm_{\inNvar}*\Poiss\left(\param,\norm_\Nvar\right)$ discussed in the Section~\ref{poisslaw}.
From its \CHF\ we see that, with $\freq=\param/\tscale$ and $D=\inNvar^2/2\tscale$, the transition
law $\norm(2Dt)*\Poiss\left(\freq t,\norm_\Nvar\right)$
\begin{equation}\label{transnp}
    \trpdf(x,t)=e^{-\freq t}\sum_{k=0}^{\infty}\frac{(\freq t)^k}{k!}\,\frac{e^{-x^2/2(k\Nvar^2+2Dt)}}{\sqrt{2\pi(k\Nvar^2+2Dt)}},
    \qquad\trchf(u,t)=e^{\freq t(e^{-\Nvar^2u^2/2}-1)}e^{-Dtu^2}
\end{equation}
The corresponding Wiener--Poisson process will have sample paths which are Brownian trajectories
interspersed with Gaussian jumps at Poisson times with intensity $\param$. The process \PDF's then
have an elementary form as time dependent Poisson mixtures of time dependent normal laws and the
corresponding process equation\refeq{FPeq} will become
\begin{equation*}
    \partial_t\prpdf(x,t)=D\partial_x^2\prpdf(x,t)+\freq\totint{[\prpdf(x+y,t)-\prpdf(x,t)]\frac{e^{-y^2/2\Nvar^2}}{\sqrt{2\pi\Nvar^2}}}{y}.
\end{equation*}
The \LS\ propagator $\norm\left(2iDt\right)*\Poiss\left(i\freq t,\norm_\Nvar\right)$ now is
\begin{equation}\label{propnp}
    \prop(x,t)=e^{-i\freq t}\sum_{k=0}^{\infty}\frac{(i\freq t)^k}{k!}\,\frac{e^{-x^2/2(k\Nvar^2+2iDt)}}{\sqrt{2\pi(k\Nvar^2+2iDt)}},\quad
    \propchf(u,t)=e^{i\freq t(e^{-\Nvar^2u^2/2}-1)}e^{-iDtu^2}
\end{equation}
and the \LS\ equation\refeq{LSeq} becomes
\begin{equation*}
    i\partial_t\wf(x,t)=-D\partial_x^2\wf(x,t)-\freq\totint{[\wf(x+y,t)-\wf(x,t)]\frac{e^{-y^2/2\Nvar^2}}{\sqrt{2\pi\Nvar^2}}}{y}.
\end{equation*}
As a second example take the law $\norm_{\inNvar}*\Poiss\left(\param,\degen_\scale\right)$
discussed in the Section~\ref{poisslaw}: from its \LCH\ $\lch(u,t)=\freq t(\cos \scale u-1)-Dtu^2$
we see that the law of the corresponding L\évy process is $\norm(2Dt)*\Poiss(\freq
t,\degen_\scale)$ and hence
\begin{eqnarray}
  \trpdf(x,t) &=& e^{-\freq t}\sum_{k=0}^{\infty}\frac{(\freq t)^k}{k!}\,\frac{1}{ 2^k}\sum_{j=0}^k\binom{k}{j}
                                \frac{e^{-[x-(k-2j)\,\scale]^2/4Dt}}{\sqrt{4\pi Dt}}\label{transdp} \\
  \trchf(u,t) &=&e^{\,\freq t(\cos\scale u-1)-Dtu^2}\nonumber
\end{eqnarray}
while the process equation\refeq{FPeq} is
\begin{equation*}
    \partial_t\prpdf(x,t)=D\partial_x^2\prpdf(x,t)+\freq\frac{\prpdf(x+\scale,t)-2\prpdf(x,t)+\prpdf(x-\scale)}{2}.
\end{equation*}
This process will have sample paths which are again Brownian trajectories interspersed with jumps
$\pm\scale$ at Poisson times with intensity $\param$. The \LS\ propagator
$\norm(2iDt)*\Poiss(i\freq t,\degen_\scale)$ instead is
\begin{eqnarray}
  \prop(x,t) &=& e^{-i\freq t}\sum_{k=0}^{\infty}\frac{(i\freq t)^k}{k!}\,\frac{1}{ 2^k}\sum_{j=0}^k\binom{k}{j}
                                \frac{e^{-[x-(k-2j)\,\scale]^2/4iDt}}{\sqrt{4\pi iDt}}\label{propdp} \\
  \propchf(u,t) &=&e^{\,i\freq t(\cos\scale u-1)-iDtu^2}\nonumber
\end{eqnarray}
and the \LS\ equation\refeq{LSeq} is
\begin{equation*}
    i\partial_t\wf(x,t)=-D\partial_x^2\wf(x,t)-\freq\frac{\wf(x+\scale,t)-2\wf(x,t)+\wf(x-\scale)}{2}.
\end{equation*}

\subsection{\Rqm\ $\relat_\scale(\freq t)$}\label{rqmequations}

We immediately see from the \CHF\ of $\relat_\scale(\param)$ that the corresponding L\évy process
$\relat_\scale(\freq t)$ will have as transition law
\begin{equation}\label{rqmtr}
    \trpdf(x,t)=
           \frac{\freq t e^{\freq
           t}K_1\left(\sqrt{\freq^2t^2+x^2/\scale^2}\right)}{\pi\scale\sqrt{\freq^2t^2+x^2/\scale^2}},
    \qquad\quad
    \trchf(u,t)=e^{\freq t\left(1-\sqrt{1+\scale^2u^2}\right)}
\end{equation}
with $\freq=\param/\tscale$ as usual. We can also explicitly write the process
equation\refeq{FPeq}
\begin{equation*}
    \partial_t\prpdf(x,t)=\freq\Zeroint{y}{[\prpdf(x+y,t)-\prpdf(x,t)]}{\frac{K_1(|y|/\scale)}{\pi|y|}}.
\end{equation*}
On the other hand the \LS\ propagator $\relat_\scale(i\freq t)$ will be given by
\begin{equation}\label{rqmprop}
    \prop(x,t)=
           \frac{i\freq t e^{i\freq
           t}K_1\left(\sqrt{-\freq^2t^2+x^2/\scale^2}\right)}{\pi\scale\sqrt{-\freq^2t^2+x^2/\scale^2}},
    \quad\qquad
    \propchf(u,t)=e^{i\freq t\left(1-\sqrt{1+\scale^2u^2}\right)}
\end{equation}
with singularities in $x=\pm\scale\freq t$, and corresponds to the \LS\ equation
\begin{equation*}
    i\partial_t\wf(x,t)=-\freq\Zeroint{y}{[\wf(x+y,t)-\wf(x,t)]}{\frac{K_1(|y|/\scale)}{\pi|y|}}.
\end{equation*}
We remember here, as remarked in the Section~\ref{relqmLaw}, that this essentially is the
integro--differential form of the well known relativistic, free Schr\"odinger equation
\begin{equation*}
    i\hbar\partial_t\wf(x,t)=\sqrt{m^2c^4-c^2\hbar^2\partial_x^{\,2}}\,\,\wf(x,t)
\end{equation*}
that we recover by taking $\freq=\param/\tscale=mc^2/\hbar$, $a=\hbar/mc$, and by reabsorbing an
irrelevant constant term $mc^2$ in a phase factor of the \WF~\cite{cufaro09}.

\section{Second kind Beta laws}\label{beta}

If $Z$ is a \rv\ with a (dimensionless) \emph{Beta law}
$\mathfrak{B}(\alpha,\beta)$ ($\alpha,\beta>0$) namely with \PDF
\begin{equation*}
    \pdf_Z(z)=\frac{z^{\alpha-1}(1-z)^{\beta-1}}{B(\alpha,\beta)}\,,\qquad0\leq
    z\leq1
\end{equation*}
then $Y=Z/(1-Z)$ is distributed according to a \emph{second kind
Beta law} $\widetilde{\mathfrak{B}}(\alpha,\beta)$ with \PDF\ (see
for example~\cite{balakrishnan})
\begin{equation*}
    \pdf_Y(y)=\frac{1}{B(\alpha,\beta)}\,\frac{y^{\alpha-1}}{(1+y)^{\alpha+\beta}}\,,\qquad0\leq
    y.
\end{equation*}
We could also introduce a scale parameter $\scale$ to get the types
$\mathfrak{B}_\scale(\alpha,\beta)$ and $\widetilde{\mathfrak{B}}_\scale(\alpha,\beta)$, but to
simplify the notation we will first consider only the dimensionless laws. Take now a third \rv\
$X=\epsilon \sqrt{Y}$ where $\sqrt{Y}$ is the \emph{positive} square root of $Y$, while $\epsilon$
is another independent \rv\ taking the two values $\pm1$ with the same probability $1/2$. We find
then that its \PDF\ is
\begin{equation*}
    \pdf_X(x)=\frac{1}{B(\alpha,\beta)}\,\frac{(x^2)^{\alpha-\frac{1}{2}}}{(1+x^2)^{\alpha+\beta}}\,.
\end{equation*}
We will use for these laws the symbol
$\widetilde{\mathfrak{B}}^{1/2}(\alpha,\beta)$ because $X$ is the
square root of a second kind Beta \rv. In particular we recover the
family of the Student laws as
\begin{equation*}
    \widetilde{\mathfrak{B}}^{1/2}\left(\frac{1}{2},\frac{\param}{2}\right)=\stud(\param),\qquad\param>0,
\end{equation*}
while $\widetilde{\mathfrak{B}}^{1/2}(3/2\,,1/2)$ is the law
introduced in the Section~\ref{cauchyevolutions} to describe the
evolution of an initial Student law $\stud(3)$ by a Cauchy
transition \PDF. For this law we have
\begin{equation*}
    \pdf(x)=\frac{2}{\pi}\,\frac{x^2}{(1+x^2)^2}
    ,\qquad\quad\chf(u)=(1-|u|)\,e^{-|u|},
\end{equation*}
and (as for the usual Cauchy laws) we find that it has neither an
expectation, nor a finite variance. The decomposition\refeq{betamix}
could now be written also as a relation within the (dimensional)
family $\widetilde{\mathfrak{B}}_\scale^{1/2}(\alpha,\beta)$ by
remembering that $\cauchy_{\vel t}=\widetilde{\mathfrak{B}}_{\vel
t}^{1/2}(1/2\,,1/2)$ and
$\stud_\inscale=\widetilde{\mathfrak{B}}_\inscale^{1/2}(1/2\,,3/2)$.
In fact, for given arbitrary scale parameters $\scale$ and
$\inscale$, and with
\begin{eqnarray*}
  P &=& \frac{1}{2}\,\frac{\scale}{\scale+\inscale} \\
  Q &=&
  \frac{1}{2}\,\frac{\scale+2\inscale}{\scale+\inscale}=\frac{1}{2}\left(1+\frac{\inscale}{\scale+\inscale}\right),
\end{eqnarray*}
we easily see that\refeq{betamix} is a special case (for $\scale=\vel t$) of
\begin{equation*}
    \widetilde{\mathfrak{B}}_\scale^{1/2}\left(\frac{1}{2},\frac{1}{2}\right)*\widetilde{\mathfrak{B}}_\inscale^{1/2}\left(\frac{1}{2},\frac{3}{2}\right)
    = P\,\widetilde{\mathfrak{B}}_{\scale+\inscale}^{1/2}\left(\frac{3}{2},\frac{1}{2}\right)
               +Q\,\widetilde{\mathfrak{B}}_{\scale+\inscale}^{1/2}\left(\frac{1}{2},\frac{3}{2}\right)
\end{equation*}
namely
\begin{eqnarray*}
    (1+\inscale|u|)e^{-(\scale+\inscale)|u|}
   &=& P\,[1-(\scale+\inscale)|u|]e^{-(\scale+\inscale)|u|}+Q\,[1+(\scale+\inscale)|u|]e^{-(\scale+\inscale)|u|} \\
  \frac{(\scale+\inscale)^2(\scale+2\inscale)+\scale x^2}{\pi[(\scale+\inscale)^2+x^2]^2}
   &=& \frac{2P}{\pi}\,\frac{(\scale+\inscale)x^2}{[(\scale+\inscale)^2+x^2]^2}
           +\frac{2Q}{\pi}\,\frac{(\scale+\inscale)^3}{[(\scale+\inscale)^2+x^2]^2}
\end{eqnarray*}

\end{document}